\documentclass[11pt]{amsart}
\usepackage{latexsym}
\usepackage{hyperref}

 \numberwithin{equation}{section}
                \usepackage{amscd}
                \usepackage{amssymb}
                \usepackage{amsmath}
                \usepackage{amsthm}
                \usepackage{enumerate}
                \usepackage{verbatim}

\def\co{\colon\thinspace}
\def\d{\delta}
\def\e{\epsilon}

\def\b{\beta}

\def\la{\lambda}

\def\tangle{\tilde{\angle}}

\def\phi{\varphi}

\def\a{\alpha}%
\def\"#1{{\accent"7F #1\penalty10000\hskip 0pt plus 0pt}} 
\def\ge{\geqslant}%
\def\le{\leqslant}%

\font\nnn=csfi10
\def\Alex{\text{{ \nnn Alex}}}

\def\Alkdv{\Alex^{n}({D, \kappa,v})}
\def\Alkd{\Alex^{n}({D, \kappa})}

\def\Alk{\Alex^{n}({ \kappa})}

\def\Z{{\Bbb Z}}
\def\R{{\Bbb R}}

\def\th{\tilde{h}}
\def\d{\delta}
\def\g{\gamma}

\newtheorem{thm}{Theorem}[section]

\newtheorem{cor}[thm]{Corollary}

\newtheorem{main}[thm]{Stability Theorem}

\newtheorem{quest}[thm]{Question}
\newtheorem{lem}[thm]{Lemma}
\newtheorem{prop}[thm]{Proposition}
\newtheorem{defn}[thm]{Definition}

\newtheorem{klem}[thm]{Key Lemma}
\newtheorem{lslem}[thm]{Local Stability Lemma}
\newtheorem{slem}[thm]{Sublemma}

\newtheorem{gluingthm}[thm]{Gluing Thorem}
\newtheorem{relgluingthm}[thm]{Relative Gluing Thorem}

\newtheorem{Example}[thm]{Example}
\newenvironment{ex}{\begin{Example}\rm}{\end{Example}}
\newtheorem{Counterexample}[thm]{Counterexample}

\newtheorem{remark}[thm]{Remark}
\newenvironment{rmk}{\begin{remark}\rm}{\end{remark}}
\newtheorem{Fact}[thm]{Fact}

\newtheorem{Nothing}[thm]{$\!\!\!$}

\newcommand{\be}{\begin{equation} }

\newcommand{\ene}{\end{equation} }

\newcommand{\ba}{\begin{eqnarray}}
\newcommand{\ea}{\end{eqnarray}}
\newcommand{\ban}{\begin{eqnarray*}}
\newcommand{\ean}{\end{eqnarray*}}

\newcommand{\vol}{\mbox{\rm vol}}

\newcommand{\curv}{\mbox{\rm curv}}

\newcommand{\diam}{\mbox{\rm diam}}

\newcommand{\Id}{\mbox{\rm Id}}

\def\conv0{{\underset{i\to\infty}{\longrightarrow}}0}
\def\conv{{\underset{i\to\infty}{\longrightarrow}}}
\def\convGH{{\overset{G-H}{\underset{i\to\infty}{\longrightarrow}}}}
\begin{document}
\abovedisplayskip=6pt plus3pt minus3pt \belowdisplayskip=6pt
plus3pt minus3pt
\title[Perelman's Stability theorem]
{Perelman's Stability theorem}


\thanks{\it 2000 AMS Mathematics Subject Classification:\rm\
53C20. Keywords: nonnegative curvature, nilpotent}\rm
\thanks{\it The  author was supported in part by an  NSF grant. 
}
\author{Vitali Kapovitch}
\address{Vitali Kapovitch\\Department of Mathematics\\University of Toronto\\
Toronto, Ontario, Canada, M5S2E4}\email{vtk@math.toronto.edu}

\begin{abstract}
We give a proof of the celebrated stability theorem of Perelman stating that for a noncollapsing sequence  $X_i$ of Alexandrov spaces with $\curv \ge k$ Gromov-Hausdorff converging to a compact Alexandrov space $X$,  $X_i$ is homeomorphic to $X$ for all large $i$.
\end{abstract}

\maketitle

\section{Introduction}

A fundamental observation of Gromov says that the class of complete  $n$-dimensional   Riemannian manifolds with fixed lower curvature and upper diameter bounds is precompact in Gromov-Hausdorff topology.
 The limit points of this class are Alexandrov spaces of dimension $\le n$ with the same lower curvature and upper diameter bounds. 
 Given a sequence of manifolds $M_i$ in the above class converging to an Alexandrov space $X$ it's interesting to know what can be said about the relationship between topologies of the limit and the elements of the sequence.

The main purpose of this paper is to give a careful proof of the following theorem of Perelman which answers this question in a more general setting of convergence of n-dimensional Alexandrov spaces in the  case when the limit space has the maximal possible dimension equal to $n$.

\begin{main}\label{intro:main}
Let $X^n$ be a compact $n$-dimensional Alexandrov space of $\curv\ge \kappa$. Then there exists an $\e=\e(X)>0$ such that for any $n$-dimensional Alexandrov space $Y^n$ of $\curv\ge \kappa$ with $d_{G-H}(X,Y)<\e$, $Y$ is homeomorphic to $X$.
\end{main}

To prove the Stability Theorem,  it's clearly  enough to show that if  $X_i^n$ is  a sequence of $n$-dimensional Alexandrov spaces  with $\curv \ge \kappa$, $\diam\le D$ converging in Gromov-Hausdorff topology to an Alexandrov  space $X$ of dimension $n$, then  Hausdorff approximations $X_i\to X$ can be approximated by homeomorphisms for all large $i$. It is this statement that will be proved in the present paper.

A proof of the Stability Theorem was given in~\cite{Per}. However, that paper is very hard to read and  is not easily accessible. 
We aim to provide a comprehensive and hopefully readable reference for Perelman's result.

Perelman also claims to have a proof of the Lipschitz version of the Stability Theorem which says that the stability  homeomorphisms can be chosen to be bi-Lipschitz. However, the proof of that result has never been written down and the author has never seen it (although he very much wants to).

It is also worth pointing out  that the Stability Theorem in dimension $3$ plays a key role in the  classification of collapsing of 3-manifolds with a lower curvature bound by Shioya and Yamaguchi~\cite{ShiYam, ShiYam2} which in turn plays a role in
Perelman's  work on  the geometrization conjecture. 
 However,  as was communicated to the author by  Kleiner,  for that particular application,
if one traces through the proofs of ~\cite{ShiYam, ShiYam2} carrying along
the additional bounds arising from the Ricci flow, then one 
finds that in each instance when a 3-dimensional Alexandrov
space arises as a Gromov-Hausdorff limit of smooth manifolds, it is in fact 
smooth, and after passing to an appropriate subsequence,
the convergence will also be smooth to a large order.  
For such convergence the  stability theorem is very well known and easily follows from Cheeger-Gromov compactness.
Therefore 
Perelman's stability theorem is unnecessary for the application 
to geometrization.

The Stability Theorem immediately implies the following finiteness theorem
 due to Grove-Petersen-Wu which was originally proved using controlled homotopy theory techniques.
\begin{thm}\label{GPW}\cite{GPW} The class of
n-dimensional Riemannian manifolds  ($n\ne 3$) with  sectional curvature $\geq
k$, diameter $\leq D$ and volume $\geq v$ has only finitely many
topological (differentiable if $n\neq 3,4$) types of manifolds.
\end{thm}

The restriction $n\ne 3$  in this theorem comes from the use of controlled homotopy results which do not hold in dimension 3.
Using the Stability theorem  one gets that homeomorphism finiteness holds in all dimensions including dimension 3. Alternatively, the homeomorphism finiteness in dimension 3 follows from an earlier result of Grove-Petersen~\cite{Gr-Pet} that the above class has finitely many homotopy types and the fact that in dimension 3, a fixed homotopy type of closed manifolds contains only finitely many homeomorphism types (this is  a direct consequence of the geometrization conjecture).

Let us mention that the proof of Stability Theorem~\ref{intro:main} presented here is fundamentally  the same as the one given in ~\cite{Per}. However, as was pointed out by Perelman in~\cite{Per-Morse}, the proof  can be simplified using the constructions developed in ~\cite{Per-Morse} and~\cite{Per-DC}. We carry out these simplifications in the present paper.

In~\cite{Per-DC} Perelman introduced the following notion:
 
 A function $f\co\R^k\to \R$ is called DC if it can be locally represented as a difference of two concave  (or, equivalently, semiconcave) functions. 
 It is easy to see that DC is an algebra and $\frac{f}{g}$ is $DC$ if both $f$ and $g$ are DC and $g$ is never zero. A map  $F\co\R^k\to\R^l$ is called DC if it has DC coordinates.   It is also easy too see (see~\cite{Per-DC}) that if $F\co\R^k\to\R^l$ and $G\co\R^l\to\R^m$ are DC then so is $G\circ F$. This allows for an obvious definition of a DC geometric structure on a topological manifold.  Despite the ease and naturality of its definition, to the best of the author's knowledge, this type of geometric structure has never been studied. In particular, the relationship between DC structures and classical geometric structures such as TOP, PL,  smooth or Lipschitz is not at all understood. It is easy to see that a PL-manifold is DC  and a DC manifold is Lipschitz but that's all one can say at the moment.
 
 The notion of  DC functions and DC homeomorphisms also makes sense on Alexandrov spaces since 
 Alexandrov spaces admit an abundance of semiconcave  (and hence DC) functions coming from distance functions. In~\cite{Per-DC}, Perelman showed that the set of regular points of an Alexandrov space has a natural structure of a DC manifold.
 
 
 The above discussion naturally leads to the following question.
 
 \begin{quest}
 Let $X_i^n$ be a noncollapsing sequence of Alexandrov spaces with $\curv\ge\kappa$, $ \diam\le D$ Gromov-Hausdorff converging to an Alexandrov space $X^n$. Is it true that $X_i$ is DC-homeomorphic to $X$ for all large $i$?

 Or more weakly.
Suppose $M_i^n$ is a noncollapsing sequence of Riemannian manifolds with $\sec\ge\kappa$, $ \diam\le D$ Gromov-Hausdorff converging to an Alexandrov space $X^n$. 
 
  Is it true that $M_i$ are  DC homeomorphic to each other for all large $i$?
 \end{quest}

Let us say a few words about the proof of the stability theorem.  
One of the main ingredients is the Morse theory for functions on Alexandrov spaces. The starting point is based on the following simple observation:

Given $k+1$ nonzero vectors in $\R^n$ with pairwise angles  bigger than $\pi/2$, any $k$ of them are linearly independent.

Motivated by this,  we'll say that   a map  $f=(f_1,\ldots,f_k)$ (with coordinates given by distance  functions $f_i=d(\cdot, a_i)$)  from an Alexandrov space to $\R^k$
 is regular at a point $p$ if there exists a point $a_{0}$ such that  the comparison angles at $p$ for the triangles $\Delta a_ipa_j$
  are bigger than $\pi/2$ for all $0\le i\ne j\le k$. 

The fibration theorem ~\cite{Per, Per-Morse} shows  that just like for regular points of smooth functions on differentiable manifolds, a map from an Alexandrov space $X$  to $\R^k$ is a topological submersion on the set of its regular points.  
In particular, if $k=\dim X$ then it's a local homeomorphism.

It's well known that for any point $p$ in an Alexandrov space,  $d(\cdot, p)$ has no critical points in a sufficiently small punctured ball $B(p,\e)\backslash \{p\}$. Therefore, by the fibration theorem, $B(p,\e)$ is homeomorphic to the cone over the metric sphere $S(p,\e)$.

The fibration theorem  (and, perhaps, even more importantly,  its proof) plays a key role in the proof of the stability theorem.  In particular it implies that Alexandrov spaces are stratified topological manifolds. 
The $k$-dimensional strata of $X$ is basically equal to the set of points $p$ in $X$ such there exist $k+1$ (but not $k+2$)  vectors in
$T_pX$ with pairwise angles bigger than $\pi/2$.

Another important tool in the proof of stability is the gluing theorem derived from deformations of homeomorphisms results of Siebenmann~\cite{Sieb}. Given a noncollapsing sequence of Alexandrov spaces $X_i\to X$,  it says that one can glue Hausdorff close  stability homeomorphisms defined on a fixed open cover of the limit space $X$ to make a global nearby homeomorphism, i.e it reduces the stability theorem to the local situation.

 To construct local stability homeomorphisms near a point $p\in X$ one argues by reverse induction on the dimension $k$ of the strata passing through $p$. The base of induction follows from the fact  mentioned above that a  map $f\co X^n\to \R^n$ is a local homeomorphism near a regular point. The general induction step is quite involved and can not be readily described in an introduction.  To show some of its flavor, we will however try to say a few words about the last step of induction from $k=1$ to $k=0$. 
 
 Given a point $p$ in the  zero strata (i.e. with the diameter of the space of directions at $p$  at most  $\pi/2$), look at a ball $B(p,r)$ where again $r$ is so small that  $d(\cdot, p)$ is regular in $B(p,r)\backslash \{p\}$.  One can show  that the same holds true for appropriately chosen lifts $p_i\in X_i$ of $p$ for all large $i$ provided $r$ is sufficiently small. This means that both $B(p,r)\backslash \{p\}$ and
  $B(p_i,r)\backslash \{p_i\}$ consist of points lying in the union of  strata of dimensions $\ge 1$.
 
 Fixing a small $\d\ll r$, the induction assumption implies the existence of homeomorphisms of the annuli $\bar{B}(p,r)\backslash B(p,\d)$ and 
 $\bar{B}(p_i,r)\backslash B(p_i,\d)$ close to the original Hausdorff approximations.  In particular, metric spheres $S(p,\d)$ and $S(p_i,\d)$ are homeomorphic.
 On the other hand, the fibration theorem implies that $B(p,\d)$ is homeomorphic to the cone over $S(p,\d)$ and $B(p_i,\d)$ is homeomorphic to the cone over $S(p_i,\d)$ for all large $i$. Gluing these homeomorphisms we obtain  homeomorphisms $B(p,r)$ onto $B(p_i,r)$ which will be close to the original Hausdorff approximations since $\d\ll r$.
 
 The general induction step is a rather nontrivial fibered  version of the above argument. 
 It is carried out in Local Stability Lemma~\ref{keylemma}.  with the main geometric ingredient provided by  Lemma~\ref{incompl}.
For technical reasons,   one has to work with more general semiconcave functions than just distance functions. An important role here is played by a technical construction from~\cite{Per-Morse} of  strictly concave functions obtained by manipulating distance functions.

Let us briefly describe  the structure of this paper.

In section \ref{sec:simple} we give  a simplified proof of the stability theorem in the special case of limits of closed Riemannian manifolds using controlled homotopy theory techniques.
In section~\ref{sec:background} we give the  background on necessary topological results on stratified spaces and some geometric constructions on Alexandrov spaces. In sections ~\ref{sec:admfun} and~\ref{s:admis} we define admissible maps and their regular points and show that they satisfy similar properties to regular points of smooth maps between manifolds.  Section~\ref{sec:proofstability} contains the proof of the stability theorem. In section~\ref{sec:subm}  we use the Stability Theorem to prove a finiteness theorem for Riemannian submersions.
In section~\ref{sec:extrstab} we generalize the Stability Theorem to show that stability homeomorphisms can be chosen to respect stratification of Alexandrov spaces into extremal subsets.

Throughout this paper we will assume some knowledge of Alexandrov geometry (see ~\cite{BGP} as a basic reference). For a  more recent treatment we strongly recommend   ~\cite{Petrunin-survey}  in the same volume.  We will also rely on the Morse theory  results and constructions from~\cite{Per-Morse} in particular the fibration theorem.
  Perelman later   generalized these results to a larger class of Morse functions in the unpublished preprint  ~\cite{Per-DC} (cf. also~\cite{Petrunin-survey}). However, these generalizations are not needed for the proof of the Stability theorem presented here.  Various references to ~\cite{Per-DC} throughout this paper are only given for bibliographical  completeness.

\bf{Acknowledgements:}\rm\  The author is profoundly grateful to Anton Petrunin and Fred Wilhelm for numerous conversations and suggestions regarding the preparation of this paper. The author would like to thank  Alexander Lytchak for  discussions and suggestions regarding the controlled homotopy theory proof of the stability theorem for manifolds.
\section{Notations}

For an open subset $U$ in a space $X$ and a subset $A\subset U$ we'll write $A\Subset U$ if $\bar{A}\subset U$.

We'll denote by   $\Alkdv$  the class of $n$-dimensional Alexandrov spaces of $\curv\ge \kappa, \diam\le D, \vol\ge v$. Similarly, we'll denote by   $\Alkd$  the class of $n$-dimensional Alexandrov spaces of $\curv\ge \kappa, \diam\le D$
and by $\Alk$ the class of $n$-dimensional Alexandrov spaces of $\curv\ge \kappa$.


For a point $v=(v_1,\ldots,v_k)\in \R^k$ and $r>0$ we'll denote by $I^k_r(v)$ the cube $\{x\in \R^k \quad |\quad |x_i-v_i|\le r \text { for all } i\}$.

For an Alexandrov space $\Sigma$ of $\curv\ge 1$ we'll often refer to its points as vectors and to distances between its points as angles.

For a metric space $X$ with $\diam\le \pi$, we'll denote by $CX$ the Euclidean cone over $X$.  Let $o$ be the vertex of $CX$. For $u\in CX$ we'll denote $|u|=d(o,u)$.

We'll call a function $h\co CX\to \R$ 1-homogeneous if $h(t\cdot x)=t\cdot h(x)$ for any $t\in \R, x\in X$.

For any space $X$ we'll denote by $KX$ the open cone on $X$ and by $\bar{K}X$ the closed cone on $X$ (i.e. the join of $X$ and a point).

For two points $p,q$ in an Alexandrov space, we'll denote by $\uparrow_p^q$ an element of $\Sigma_pX$ tangent to a shortest geodesic connecting $p$ to $q$. We'll denote by $\Uparrow_p^q$ the set of {\it all} such directions.

For three points $x, p, y$ in an Alexandrov space $X$ of $\curv\ge\kappa$,  we'll denote by $\tilde{\angle} xpy$ the comparison angle at $p$, i.e the angle $\angle\tilde{x}\tilde{p}\tilde{y}$ in the triangle $\tilde{x}\tilde{p}\tilde{y}$ in the complete simply connected space of constant curvature $\kappa$ with $d(x,y)=d(\tilde{x},\tilde{y})$, $d(p,y)=d(\tilde{p},\tilde{y})$, $d(x,p)=d(\tilde{x},\tilde{p})$, 

We will also often use the following convention.
In the proofs of various theorems we'll denote by $c$ or $C$  various constants  depending on the dimension and the lower curvature bound present and which sometimes will depend on additional parameters present. When important this dependence will be clearly indicated.  

We will denote by $\varkappa$ various continuous  increasing   functions  $\varkappa\co R_+\to\R_+$ satisfying  $\varkappa(0)=0$.

By $o(i)$ we will denote various  positive functions on $\Z_+$ such that $o(i)\conv0$.

We'll write  $o(i|c)$ to indicate a function which depends on an extra parameter $c$ and satisfies  $o(i|c)\conv0$ for any fixed $c$.
Sometimes we'll use the same convention for $\varkappa(\d|c)$.

\section{Simplified proof of the Stability Theorem for limits of manifolds}\label{sec:simple}

The author is grateful to A.~Lytchak for bringing to his attention the fact,  that using todays knowledge of the local structure of Alexandrov spaces, a simple proof of the stability theorem  for $n\ge 4$ can be given for the special case of  limits of Riemannian manifolds.

The proof uses controlled homotopy theory techniques employed in~\cite{GPW}.

Let us briefly describe the argument the general otline of which was suggested to the author by A.~Lytchak. 

It is now well-known that the class  $\Alkdv $   has a common contractibility function (see e.g. ~\cite{Per-Pet} or \cite{Per-Morse}). Therefore, by~\cite[Lemma 1.3]{GPW},  if 
 $X^n_i\convGH X^n$ is a convergent sequence  in  $\Alkdv $ , then $X_i^n$ is  $o(i)$-homotopy equivelent to $X$ for all large $i$ (see  \cite{GPW}). This means that there are  homotopy equivalences $f_i\co X_i\to X$ with  homotopy inverses $h_i\co X\to X_i$ such that $f_i\circ h_i\simeq \Id_X$ and $h_i\circ f_i\simeq \Id_{X_i}$ through homotopies $F_i\co X\times[0,1]\to X $ and 
 $G_i\co X_i\times[0,1]\to X_i$  
 such that all the point trajectories of $F_i$ and $f_i\circ G_i$ have 
 $o(i)$-small diameters in $X$. In fact, in our case, one can make $f_i, h_i$ to be $o(i)$-Hausdorff approximations.

It is relatively easy to show (see~\cite{GPW}) that  if  $M_i\in  \Alkdv $  is a sequence of closed  Riemannian manifolds Gromov-Hausdorff converging to a space $X$,  then $X$ is a homology manifold.

 However, at the time of the writing of ~\cite{GPW} the local structure of Alexandrov spaces was not well understood and it was therefore not known if  $X$ is an actual manifold. This made the application of controlled topology techniques employed in ~\cite{GPW} fairly tricky. 

We first show that $X$  is a manifold.

We'll need the following result from~\cite{Kap} (cf. \cite{GrWi}).
\begin{lem}\label{convex-lift}
Suppose  $M_i\in   \Alkdv $  is a sequence of  Riemannian manifolds Gromov-Hausdorff converging to a space $X$.
Let $p\in X$ be any point. 

 Then there exists  a $\d>0$,   a 1-Lipschitz function $h\co X\to\R$, strictly convex on $B(p,\d)$  such that $h(p)=0$ is a strict local minimum of $h$ and a sequence of smooth  1-Lipschitz  functions $h_i\co M_i\to\R$ uniformly converging to $h$ such that $h_i$ is
 strictly convex on $B(p_i,\d)$ (where $p_i\in M_i$ converges to $p$)  for all large $i$.
\end{lem}

This lemma together with the fibration theorem of Perelman~\cite{Per-Morse} easily implies that a noncollapsing  limit of Riemannian manifolds with lower sectional curvature bound is  a topological  manifold.

\begin{lem}\label{man-lem}
Let $M_i\in  \Alkdv $ be a sequence of closed Riemannian manifolds converging to an Alexandrov space $X$.
Then $X$ is a topological manifold.
\end{lem}
\begin{proof}
We argue by induction on dimension.  
The cases of $n\le 2$ are easy and are left to the reader as an exercise.

Suppose $n\ge 3$.  Let $p\in X$ be any point and let $h$ be the function provided by Lemma~\ref{convex-lift}.
Then for some small $\e>0$ the set $\{h\le \e\}$ is a compact convex subset of $X$ and $\{h_i\le\e\}$ is compact convex in $M_i$ for large $i$. Obviously, $\{h_i\le\e\}\convGH \{h\le \e\}$. By~\cite[Theorem 1.2]{Ptr1}, we have that  $\{h_i=\e\}\convGH \{h= \e\}$ with respect to their induced inner metrics. By the Gauss Formula, $\sec(\{h_i=\e\})\ge\kappa$  for all large $i$ and hence $ \{h= \e\}$ is an  $(n-1)$-dimensional  Alexandrov space of $\curv\ge\kappa$ which is a manifold by induction assumption. Moreover, since $h_i$ is smooth, strictly convex with unique minimum, $ \{h_i= \e\}$ is diffeomorphic to $S^{n-1}$ for large $i$ and hence $ \{h= \e\}$ is a homotopy $(n-1)$-sphere.
 Since $h$ is strictly convex in $\{h\le \e\}$,  it has no critical points in $\{h\le \e\}\backslash \{p\}$. Therefore, by~\cite[Theorem 1.4]{Per-Morse}, $\{h<\e\}$ is homeomorphic to the open cone over $\{h= \e\}$.
 
 If $n=3$ then $ \{h= \e\}$ is obviously homeomorphic to $S^2$ which means that $\{h<\e\}$ is homeomorphic to $\R^3$.
 If $n=4$ then $ \{h= \e\}$ is a homotopy $3$-sphere and a manifold. By the work of Freedman~\cite[Corollary 1.3]{Freed}, this implies that the cone over $ \{h= \e\}$ is homeomorphic to $\R^4$. 
 If $n\ge 5$ then $ \{h= \e\}$ is a homotopy $S^{n-1}$ and a manifold and hence is homeomorphic to $S^{n-1}$ by the Poincare conjecture.

  \end{proof}
  
  \begin{rmk}
  As was suggested to the author by A.~Lytchak, a different proof of  Lemma \ref{man-lem} can be given for $n\ge 5$ by verifying that $X$ satisfies the disjoint disk  property and thus is a manifold by a result of Edwards~\cite{Dav}. However, the author prefers his own argument given above.
\end{rmk}
Now that we know that $X$ is a manifold and $f_i\co M_i\to X$ are $o(i)$ -homotopy equivalences, we can apply the results from   ~\cite{Chap-Fer} for $n\ge 5$ and ~\cite{QuinnIII} together with ~\cite{Chap-Fer}  for $n=4$, which say that  under these conditions $f_i$s can be $o(i)$- approximated by  homeomorphisms for all large $i$.

The same holds in dimension 3 by~\cite{Jak} but only modulo the Poincare conjecture.

However,  it is possible  (and would certainly be  a lot more preferable) that one can  use the fact that all the subsets $\{h_i\le \e\}$ are actually topological balls and not merely contractible to give a proof in dimension 3 which does not rely on the Poincare conjecture.

\begin{rmk}
 By using relative versions of controlled homotopy theory results  from  ~\cite{QuinnIII} and  ~\cite{Chap-Fer} mentioned above, it should be possible to generalize the above proof to the case of pointed Gromov-Hausdorff  convergence of manifolds.
This would 
amount to the manifold case of Theorem~\ref{pointedstab} below.
 Alternatively, one can handle the pointed case as follows.

Suppose we have a pointed convergence $(M_i^n,q_n)\to (X^n,q)$ where $M^n_i$ are (possibly noncompact)  Riemannian manifolds of $\sec\ge \kappa$. The proof of Lemma~\ref{man-lem} is obviously local and hence $X$ is a topological manifold. 

Let $p$ be any point in $X$. 
Let $h$ and $h_i$ be the functions constructed in the proof of Lemma~\ref{man-lem}. Let $Y_i$ be the double of $\{h_i\le \e\}$ and $Y$ be the double of $\{h\le \e\}$.  Obviously $Y_i\convGH Y$ and $Y_i$ is homeomorphic to $S^n$ for large $i$.  While the metric on $Y_i$ is not smooth along the boundary of  $\{h_i\le \e\}$, it's easy to see that the proof of Lemma~\ref{man-lem} still works for the convergence $Y_i\to Y$ and hence $Y$ is a closed   topological manifold.  
By the same controlled homotopy theory results used earlier, we conclude that Hausdorff approximations  $Y\to Y_i$ can be $o(i)$-approximated by homeomorphisms $g_i\co Y\to Y_i$ for all large $i$. Restricting $g_i$ to  $\{h< \e/2\}$ we obtain an open embedding of  $\{h< \e/2\}$ into $M_i$, which is $o(i)$-close to the original Hausdorff approximation $(M,q)\to (M_i,q_i)$.

Finally, by using  topological gluing theorem~\ref{locgluing}  below, for any fixed $R>0$ we can glue finitely many such local homeomorphisms to get an open embedding $B(q,R)\to Y_i$, which is $o(i)$-Hausdorff close to the original Hausdorff approximation  $(M,q)\to (M_i,q_i)$.

\end{rmk}

\section{Background}\label{sec:background}
\subsection{Stratified spaces} 
Most of the material of this section is taken with almost no changes  from ~\cite{Per} as no significant simplifications or improvements of the exposition seem to be possible.
\begin{defn}
A metrizable space $X$ is called an MSC-space (space with multiple conic singularities) of dimension $n$ if  every point $x\in X$ has a neighborhood pointed homeomorphic to an open cone over a compact  $(n-1)$-dimensional MCS space. Here we assume the empty set to be the unique $(-1)$-dimensional MCS-space.
\end{defn}

\begin{rmk}
A compact $0$-dimensional MCS-space is a finite collection of points with discrete topology.
\end{rmk}

\begin{rmk}
An open conical neighborhood of a point in an MCS-space is unique up to pointed homeomorphism ~\cite{Kwun}.
\end{rmk}

It easily follows from the definition that an MCS space  has a natural topological stratification.

We say that a point $p\in X$ belongs to the $l$-dimensional strata $X_l$ if $l$ is the maximal number $m$ such that the conical neighbourhood 
of $p$ is pointed homeomorphic to $R^m\times K(S)$  for some  MCS-space $S$. It is clear that $X_l$ is an $l$-dimensional topological manifold.

We will need two general topological results which hold for spaces more general than  Alexandrov spaces and follow from the general theory of deformations of homeomorphisms developed by Siebenmann~\cite{Sieb}.

\begin{thm}\label{local bundle}\cite[Theorem 5.4, Corollary 6.14, 6.9]{Sieb}

Let $X$ be a metric space and $f\co X\to \R^k$ be a continuous, open, proper map such that for each $x\in X$ we have

\begin{enumerate}
\item $f^{-1}(f(x))$ is a compact MCS-space;
\item $x$ admits a product neighborhood with respect to $f$, i.e there exists an open neighbourhood $U_x$ of $x$ and a homeomorphism $F_x\co U_x\to U_x\cap f^{-1}(f(x))\times f(U_x)$ such that $f_x=p_2\circ F_x$ where $p_2\co U_x\cap f^{-1}(f(x))\times f(U_x)\to \R^k$ is the coordinate  projection onto the second factor.
\end{enumerate}

Then $f$ is a locally trivial fiber bundle.

Moreover, suppose we have  in addition   that $f(U_x)=I^k$. Let $K\subset U_x$  be a compact subset.  Then there exists a homeomorphism
$\phi\co f^{-1}(I^k)\to f^{-1}(f(x))\times I^k$ respecting $f$ (I.e. such that $f=p_2\circ \phi$. and such that $\phi|_K=F_x|_K$.
\end{thm}

The next gluing theorem is the key topological ingredient in the proof of the Stability Theorem.   It says that for MCS spaces close local homeomorphisms given on a finite  open covering can be glued to a nearby  global homeomorphism under some mild (but important!) geometric assumptions.

First we need a technical definition.

\begin{defn}
A metric space $X$ is called $\varkappa$-connected if for any two points $x_1,x_2\in X$ there exists a curve  connecting $x_1$ and $x_2$ of $\diam\le \varkappa(d(x_1,x_2))$.
\end{defn}

\begin{gluingthm}\label{gluingthm}
Let $X$ be a compact MCS-space, $\{U_\a\}_{\a\in \frak A}$ be a finite covering of $X$. Given a function $\varkappa_0$, there exists $\varkappa=\varkappa(X, \{U_\a\}_{\a\in \frak A},\varkappa_0)$ such that the following holds:

Given a $\varkappa_0$-connected MCS-space $\tilde{X}$,  an open cover of $\tilde{X}$  $\{\tilde{U}_\a\}_{\a\in \frak A}$, a $\delta$-Hausdorff approximation $\phi\co X\to \tilde{X}$ and a family of homeomorphisms $\phi_\a\co  U_\a \to \tilde{U}_\a$, $\delta$-close to $\phi$,

then there exists a homeomorphism $\bar{\phi}\co X\to  \tilde{X}$, $\varkappa(\d)$-close to $\phi$.
\end{gluingthm}

\begin{proof}[Proof of Theorem~\ref{gluingthm}]
This proof of Theorem~\ref{gluingthm} is taken verbatim from~\cite{Per}.

We'll need two lemmas.

\begin{lem}[Deformation Lemma]\label{deflem}
Let $X$ be a compact metric MCS-space,  $W\Subset V\Subset U\subset X$ be open subsets.

Let $\phi\co U\to X$ be an open embedding $\d$-close to the inclusion $i$.

Then there exists an open embedding $\psi\co U\to X$, $\varkappa(\d)$-close to $i$ and such that  $\psi\equiv\phi$ on $W$ and $\psi\equiv i$ on $U\backslash V$.

(Here $\varkappa$ depends on $W,U,V,X$ but not on $\phi$).
\end{lem}
\begin{proof}
Consider the open embedding $\phi |_{U\backslash \bar{W}}\co U\backslash \bar{W}\to X$. By the deformation theorem of Siebenmann~\cite[Theorem 5.4]{Sieb},  it can be perturbed to an open embedding $\phi_1\co U\backslash \bar{W}\to X$ which is $\varkappa(\d)$-close to $i$,  coincides with $i$ on some neighborhood of $\partial V$ and is equal to $\phi$ outside some compact subset of  $U\backslash \bar{W}$. Now let

\[
\psi(x)=
\begin{cases}
\phi(x),\text{ for } x\in W\\
\phi_1(x), \text{ for } x\in V\backslash W\\
x, \text{ for } x\in U\backslash V
\end{cases}
\]

It's clear that $\psi$ satisfies the conclusion of the Lemma.

\end{proof}
\begin{lem}\label{ontolem}
Under the assumptions of Theorem~\ref{gluingthm},  let $x\in X, \tilde{x}\in\tilde{X}$ satisfy $d(\phi(x),\tilde{x})<\d$. Let $V\subset X$ be an open set containing $B(x,\varkappa_0(\d)+10\d)$.
Suppose $\psi\co V\to\tilde{X}$ be an open embedding, $\d$-close to $\phi$.

Then $\tilde{x}\in \psi(V)$.

\end{lem}
\begin{proof}
Let $\tilde{\g}\co [0,1]\to \tilde{X}$ be a curve of $\diam\le \varkappa_0(\d)$ with $\tilde{\g}(0)=\psi(x),\tilde{ \g}(1)=\tilde{x}$.

We'll  show that $\tilde{\g}$  has a lift $\g\co [0,1]\to V$ with respect to $\psi$. Since $\psi $ is an open embedding we can lift $\tilde{\g}$ on some interval $[0,\e)$. Observe that given a lift $\g$  of $\tilde{\g}$  on $[0,t)$ for some $t\le 1$  it can always be extended to $[0,t]$ provided the closure of $\g([0,t))$ is contained in $B(x,\varkappa_0(\d)+10\d)$. The fact that $\psi$ is $\d$-close to $\phi$ and $\phi$ is a $\d$-Hausdorff approximation assures that this is always the case. Therefore the lift $\g$ can be extended to $[0,1]$ with $\psi(\g(1))=\tilde{x}$.
\end{proof}

\begin{rmk}
The proof of Lemma~\ref{ontolem} is the only place in the proof of Theorem~\ref{gluingthm} where we use the assumption that $\tilde{X}$ is $\varkappa_0$-connected.
\end{rmk}

We can now continue with the proof of Theorem~\ref{gluingthm}.

Suppose $U_{\a_1}\cap U_{\a_2}\ne \emptyset$. Let $U_1^4\Subset U_1^3\Subset U_1^2\Subset U_1^1\Subset U_{\a_1}$ and 
$U_2^4\Subset U_2^3\Subset U_2^2\Subset U_2^1\Subset U_{\a_2}$ be open subsets such that $U_1^4, U_2^4$ still cover
$X\backslash \cup_{\a\in \frak A\backslash\{\a_1,\a_2\}}U_\a$.

By Lemma~\ref{ontolem},  we have $\phi_{\a_1}(U_1^1\cap U_2^1)\subset \phi_{\a_2}(U_{\a_2})$ provided $\d$ is small enough.

Therefore we can consider the open embedding $\phi_{\a_2}^{-1}\circ \phi_{\a_1}\co U_1^1\cap U_2^1\to U_{\a_2}$. Clearly, it is $2\d$-close to the inclusion $i$. By Lemma~\ref{deflem}, there exists an open embedding $\psi\co U_1^1\cap U_2^1\to U_{\a_2}$, $\varkappa(\d)$-close to $i$ and such that $\psi\equiv \phi_{\a_2}^{-1}\circ \phi_{\a_1}$ on $U_1^3\cap U_2^3$ and $\psi\equiv i$ on $U_1^1\cap U_2^1\backslash U_1^2\cap U_2^2$. We can extend $\psi$ to $U_2^1$ by setting $\psi\equiv i$ on $U_2^1 \backslash  U_1^2\cap U_2^2$ and define $\phi_{\a_2}'=\phi_{\a_2}\circ\psi$.

Now we define $\phi'\co U_1^4\cup U_2^4\to\tilde{X}$ by the formula
\[
\phi'(x)=
\begin{cases}
\phi_{\a_1}(x) \text{ for } x\in U_1^4\\
\phi_{\a_2}'(x) \text{ for } x\in U_2^4
\end{cases}
\]

It's clear that $\phi'$ is an open immersion and it's is actually an embedding provided $\d$ is small enough.

Moreover, by Lemma~\ref{ontolem} we have 
$$\tilde{X}\backslash \underset{\a\in \frak A\backslash\{\a_1,\a_2\}}  { \cup}\tilde{U}_\a\subset \phi'(U_1^4\cup U_2^4).$$

 Now the statement of the theorem immediately follows by induction on the number of elements in $\frak A$.

\end{proof}
In  fact, we will need a somewhat stronger version of this theorem which assures that the gluing can be done relative to a fiber bundle structure on all the limit and approximating spaces.

\begin{thm}[Strong Gluing Theorem]\label{stronggluing}
Under the assumptions of Gluing Thorem~\ref{gluingthm} we are given in addition continuous maps $f\co X\to \R^k, \tilde f\co \tilde{X}\to \R^k, h\co X\to R, \tilde{h}\co \tilde{X}\to\R$ and a compact set $K\subset X$ such that the following holds
\begin{enumerate}
\item for any $U_\a$ with $U_\a\cap K\ne \emptyset$ we have $(\tilde{f},\tilde{h})\circ \phi_\a=(f,h)$
\item for any $U_\a$ with $U_\a\cap K=\emptyset$ we have $\tilde{f}\circ \phi_\a=f$
\item for any  $U_\a$ with $U_\a\cap K\ne \emptyset$ , $U_\a$ is contained in a product neighbourhood with respect to $(f,h)$
\item for any $U_\a$ with $U_\a\cap K=\emptyset$,  $U_\a$ is contained in a product neighbourhood with respect to $f$
\end{enumerate}

Then the gluing homeomorphism $\bar{\phi}$ can be chosen to respect $f$ on $X$ and $(f,h)$ on $K$  (i.e $(\tilde{f},\tilde{h})\circ \bar{\phi}=(f,h)$ on $K$ and  $\tilde{f}\circ \bar{\phi}=f$ on $X$.

\end{thm}

\begin{proof}

The proof of Theorem~\ref{gluingthm} can be trivially adapted to prove  Strong Gluing Theorem~\ref{stronggluing}  once we observe that  the  theorem of Siebenmann quoted in the proof of Lemma~\ref{deflem} has a stronger version respecting products with $\R^k$~\cite[Theorem 6.9]{Sieb} so that the deformation $\psi$ given by   Lemma~\ref{deflem} can be made to respect the product structure $X\cong X_1\times \R^k$ if $\phi\co U\to X$ respects that product structure.

\end{proof}

For applications to pointed Gromov-Hausdorff convergence we will need the following  local version of the Gluing theorem for which the requirement that the approximated space be $\varkappa$-connected can be slightly weakened. For simplicity, we only state the unparameterized version.

\begin{thm}\label{locgluing}

Let $U\Subset V\Subset W\subset X$ be relatively compact  open subsets in an MCS-space $X$. Let 
$\{U_\a\}_{\a\in \frak A}$ be a finite covering of $\bar{W}$ with the property that if $U_\a\cap \bar{V}\ne \emptyset$ then $U_\a\Subset W$.

Then  given a function $\varkappa_0$, there exists $\varkappa=\varkappa(X, U,V,W, \{U_\a\}_{\a\in \frak A},\varkappa_0)$ such that the following holds:

Given a $\varkappa_0$-connected MCS-space $X'$,  and  subsets $U'\Subset V'\Subset W'\Subset X'$,  an open cover  $\{U'_\a\}_{\a\in \frak A}$ of $\bar{W}'$ , a $\delta$-Hausdorff approximation $(\bar{W}, \bar{V}, \bar{U})\to (\bar{W}', \bar{V}', \bar{U}')$ and a family of homeomorphisms $\phi_\a\co  U_\a \to {U}'_\a$, $\delta$-close to $\phi$,

then there exists an open embedding  $\phi'\co V\to X'$, $\varkappa(\d)$-close to $\phi$ such that $\phi(V)\supset U'$ if $\d$ is sufficiently small.
\end{thm}

\begin{proof}
The proof is exactly the same as the proof of Theorem~\ref{gluingthm} except in the induction procedure we only glue the embeddings of those $U_\a$ which intersect $\bar{V}$.
\end{proof}

\begin{defn}
A map $f\co X\to Y$  between two metric spaces is called $\e$-co-Lipschitz if for any $p\in X$ and all small $R$ we have $f(B(p,R))\supset B(f(p),\e R)$.
\end{defn}

We will often make use of the following simple observation the proof of which is left to the reader as an exercise.

\begin{lem}\label{colip}
Let $f\co X\to Y$ is $\e$-co-Lipschitz where $X$ is compact. Let $p\in X$ and $\g\co [0,1]\to Y$ be a rectifiable curve with $\g(0)=f(p)$. Then there exists  a lift $\tilde{\g}\co [0,1]\to X$ of $\g$  such that $\tilde{\g}(0)=p$ and $L(\tilde{\g})\le \frac{1}{\e}L(\g)$.

\end{lem}

\subsection{Polar vectors}
\begin{defn}
Let $\curv \Sigma\ge 1$.  
 Given  elements $u,v\in C\Sigma$  we define $\langle u, v\rangle$ by the formula
\[
\langle u, v\rangle= |v|\cdot |u|\cdot \cos\angle uv
\]
\end{defn}

\begin{defn}
Let $\curv \Sigma\ge 1$.  
 Two vectors $u,v\in C\Sigma$ are called {\sf polar} if for any $w\in C\Sigma$ we have
\[
\langle v, w\rangle+\langle u, w\rangle\ge 0
\]

More generally,  $u$ is called polar to a set $V\subset C\Sigma$ if for any $w\in C\Sigma$ we have

\[
\sup_{v\in V}\langle v, w\rangle+\langle u, w\rangle\ge 0.
\]

\end{defn}

It is known (see \cite{PP}  or \cite[Lemma 1.3.9]{Petrunin-survey}), that for any $v\in \Sigma$ there exists $u$, polar to $v$.

A function $f\co \Sigma\to  \R$ is called {\it spherically concave} if for any $y$ lying on a shortest geodesic connecting $x$ and $z\in \Sigma$ we have

\[
d(x,z)f(y)\ge d(x,y)f(z)+d(z,y)f(x).
\]

As with ordinary concave functions, for a space with boundary we demand that the canonical extension of $f$ to the doubling of $\Sigma$ be spherically concave.

It's easy to see that $f$ is spherically concave iff its  1-homogeneous extension to $C\Sigma$ is concave.

We will  need the following property of polar vectors~\cite[Section 1.3.8]{Petrunin-survey}:

Let $f\co C\Sigma\to\R$ be  concave and 1-homogeneous. Suppose 
$u\in C\Sigma$ is  polar to $V\subset C\Sigma$. Then 
\begin{equation}\label{e:polar}
f(u)+\inf_{v\in V} f(v)\le 0.
\end{equation}

Finally, we'll make use of the following fact~\cite[Section 1.3.8]{Petrunin-survey}:

Given any two distinct points $p,q$ in an Alexandrov space  we have 
\begin{equation}\label{e:polarg}
\nabla_p d(\cdot, q) \text{ is polar to } \Uparrow_p^q.
\end{equation}
\subsection{Gradient flows of semiconcave functions}
$\qquad$

It was shown in ~\cite{PP}  (cf. ~\cite{KPT} ) that semiconcave  functions on Alexandrov spaces admit well defined forward gradient flows.

Moreover,  one can bound the Lipschitz constant of the gradient flows as follows.

Suppose $f\co X\to \R$ is $\la$-concave. Let $\alpha$ and $\beta$ are two $f$-gradient curves. 
 Then
 \begin{equation}\label{e:lip}
d(\alpha(t_1),\beta(t_1))\le d(\alpha(t_0),\beta(t_0))\exp(\lambda (t_1-t_0))
\ \ \text{for all}\ \ t_1\ge t_0
\end{equation}

In particular, if $f$ is concave then its gradient flow $\phi_t$ is $1$-Lipschitz for any $t>0$.

\section{Admissible functions and their derivatives}\label{sec:admfun}
\begin{defn}
Let $X$ be an Alexandrov space of $\curv\ge\kappa$.
Let $f\co X\to\R$ have the form $f=\sum_\a\la_\a\phi_\a(d(\cdot, A_\a))$  where each $A_\a$ is a closed subset of $X$, $\la_\a\ge 0, \sum_\a\la_\a\le 1$ and each $\phi_\a\co \R\to\R$ is a twice differentiable function with $0\le \phi_\a'\le 1$. 
We say that such $f$ is {\sf admissible}  on $U=X\backslash \cup_\a A_\a$.
\end{defn}

With a slight abuse of notations we'll sometimes simply say that $f\co X\to\R$ of the above form is admissible.

It's obvious from the definition that an admissible function is $1$- Lipschitz and semiconcave on $U$. More precisely,  $f$ is $\la$-concave near  $p\in U$ where $\la$ depends on $\kappa$ and $d(p,\cup_iA_i)$.

By the first variation formula, 
$df_p\co\Sigma_pX\to\R$  has the form
\[
df_p=\sum_\a -a_\a\cos(d(\cdot, \Uparrow_p^{A_\a})) \text{  where } a_\a=\la_a\phi'_\a(d\cdot, A_\a)).
\]

Note that $a_\a\ge 0, \quad \sum_\a a_\a\le 1$.

In  view of this, following~\cite{Per-Morse}, we call a function   $h\co \Sigma\to\R$ where $\curv\Sigma\ge 1$ a function of class DER
if it has the form

\[
\sum_\a -a_\a\cos(d(\cdot, {A_\a})) \text{  where } a_\a\ge 0, \quad \sum_\a a_\a\le 1.
\]

for some  finite collection $\{A_\a\}_{\a\in \frak A}$ of subsets of $\Sigma$.

DER functions are spherically concave on $\Sigma$ and concave when radially extended to $C\Sigma$.

We define the scalar product of two functions in DER by the formula:
for
\begin{equation}\label{e:rep}
h=\sum_\a-a_\a\cos(d(\cdot,A_\a)), \quad g=\sum_\b-b_\beta\cos(d(\cdot,B_\beta))
\end{equation}
put 

\[
\langle h, g\rangle=\sum_{\a,\b}a_\a b_\b \cos(d(A_\a,B_\b)).
\]

Note that this definition depends on the representations of $h,g$ given by (\ref{e:rep}) and not just the values of $f,g$ at every point.

It is shown in ~\cite{Per-Morse} that the scalar product satisfies the following properties:

\begin{enumerate}[(i)]
\item $ \langle d_ph,d_pg\rangle\le \langle h, g\rangle-h(p)g(p)$ for any $p\in\Sigma$;
\item $ \langle h,h \rangle\ge \left(\inf_{g\in DER(\Sigma)}\langle h, g\rangle \right)^2\ge  0$;
\item For any $h\in DER(\Sigma)$ there is a point $\hat{A}\in\Sigma, 0\le \hat{a}\le 1$ such that for
 $\hat{h}=-\hat{a}\cos(d(\cdot, \hat{A}))$ we have $\langle h, g\rangle\ge \langle \hat{h}, g\rangle$.
\end{enumerate}

Notice that an admissible function $f$ on an Alexandrov space $X$ can be naturally lilfted to a nearby Alexandrov space $\tilde{X}$ by lifting the sets $A_\a$ from the definition of $f$ to nearby sets in $\tilde{X}$ and defining $\tilde{f}$ by the same formula as $f$.

We will need the following simple lemma:

\begin{lem}\label{l:limsup}
Let $X_i\convGH X$ where $\curv X_i\ge \kappa$ and let $f,g\co X\to \R$ be  admissible at $p\in X$.  Let $f_i,g_i\co X_i\to\R$ be natural lifts of $f,g$. Let $X_i\ni  p_i\conv p$. Then
\[
\langle d_pf,d_pg\rangle \ge \limsup_i \langle d_{p_i}{f_i},d_{p_i}{g_i}\rangle.
\]
\end{lem}
\begin{proof}
By linearity it's easy to see that it's enough to prove the lemma for $f,g$ of the form $f=d(\cdot, A), g=d(\cdot, B)$. For functions of this form the statement easily follows from Toponogov comparison  by an argument by contradiction.
\end{proof}

\section{Admissible  maps and their regular points}\label{s:admis}

\begin{defn}
A map $g\co X\to\R^k$ is called admissible on an open set $U\subset X$ if it admits a representation $g=G\circ f$ where all the components of $f\co U\to \R^k$ are admissible on $U$ and $G$ is a bi-Lipschitz homeomorphism between open sets in $\R^k$.
\end{defn}
\begin{rmk}
The inclusion of the bi-Lipschitz map $G$ in the definition of an admissible map might seem strange at this point.  
However, it will  significantly simplify certain steps in the proof of the Stability theorem as well as prove useful in applications.

To unburden the exposition we will employ the following convention.
If an admissible map is denoted by $f$ (with any indices) we'll automatically assume that in the above definition $G\equiv \Id$.

\end{rmk}

\begin{defn}\label{reg-map}
An admissible map $g\co U\to \R^k$   is called $\e$-regular at $p\in U$ if for some representation $g=G\circ f$ of $g$ the following holds:

\begin{enumerate}
\item $\min_{i\ne j}-\langle d_pf_i,d_qp_j\rangle  >\e$;
\item There exists $v\in \Sigma_pX$ such that $f_i'(v)>\e$ for all $i$.
\end{enumerate}

We say that $g$ is regular at $ p$ if it's $\e$-regular at $p$ for some $\e>0$.
\end{defn}

It is easy to see that $f=(f_1,\ldots,f_k)$ is $\e$-regular at $p$ iff there is a point $q$ such that for $f_0=d(\cdot, q)$ we have

\[
\langle d_pf_i,d_pf_j\rangle  <-\e \text{ for all }i\ne j.\]

\begin{ex}
As was mentioned in the introduction, a  basic example of a regular map is as follows. Suppose $\{f_i=d(\cdot, A_i)\}_{i=0,\ldots,k}$ satisfy
\[
\angle\Uparrow_p^{A_i}\Uparrow_p^{A_j}>\pi/2 \text{ for all }i\ne j.
\]

Then $f=(f_1,\ldots,f_k)$ is regular at $p$.

This example shows that regular points of admissible maps naturally generalize regular points of smooth maps because of the following simple observation:

Given $k+1$ non-zero vectors in $\R^n$ with all pairwise angles $>\pi/2$, any $k$ of them are linearly independent.
\end{ex}

As an obvious corollary of  Lemma~\ref{l:limsup} we obtain

\begin{cor}\label{reg-lifts}
Let $f\co X\to \R^k$ be admissible and $\e$-regular at $p\in X$.   Let $0<\e'<\e$. Suppose $X_i\convGH X$ where $X_i\in \Alkd$.  Let  $X_i\ni  p_i\conv p$ and let $f_i\co X_i\to\R^k$ be  natural lifts of $f$. Then there exists $\d>0$ such that $f_i$ is $\e'$-regular on $B(p_i,\d)$ for all large $i$.
\end{cor}

In particular,  the set of $\e$-regular (regular) points of an admissible map is open.

\begin{rmk}
In ~\cite{Per-DC}, Perelman generalized all the Morse theory results from ~\cite{Per-Morse} to the more general and much more natural class of admissible fuctions consisting of 1-Lipschitz semiconcave functions.
Corollary~\ref{reg-lifts} is one of the main reasons why we restrict the class of admissible functions to the rather special semiconcave functions constructed from distance functions.  Various  definitions of regularity are possible for maps with semiconcave coordinates  (see~\cite{Per-DC,Petrunin-survey}).  While all these different definitions allow for relatively straightforward generalization of the results from ~\cite{Per-Morse}, the author was unable to prove the analogue of Corollary ~\ref{reg-lifts} using any of these definitions.

Another (perhaps more serious) reason why we are are forced to work with a  small class of admissible functions is that there is currently no known natural way
of lifting general semiconcave functions from the limit space to the elements of the sequence.
\end{rmk}

The following Lemma is due to Perelman.

\begin{lem}\cite[Lemma 2.3]{Per-Morse},~\cite[Lemma 2.2]{Per-DC}, \cite[Lemma 8.1.4]{Petrunin-survey}\label{per-lem}

Suppose $\Sigma^{n-1}$ has $\curv\ge 1$ and let $f_0,\ldots, f_k\co C\Sigma\to\R$ be 1-homogeneous concave functions such that
$\e=\min_{i\ne j}-\langle f_i,f_j\rangle  >0$.
Then

\begin{enumerate}
\item $k\le n$;
\item There exists $w\in \Sigma$ such that $f_i(w)> \e$ for all $i\ne 0$.
\item There exists $v\in \Sigma$ such that $f_0(v)>\e, f_1(v)<-\e$ and $f_i(v)=0$ for $i=2,\ldots,k$.
\end{enumerate}
\end{lem}

Let $g\co X\to R^k$ be an admissible map. Let $X^\e_{reg}(g)\subset X$  ($X_{reg}(g)\subset X$) be the set of  $\e$-regular (regular) points of $g$.

Then the following properties hold~\cite{Per-Morse, Per-DC}
\begin{enumerate}[a)]
\item $X^\e_{reg}(g)$ (and hence also ($X_{reg}(g)$)  is open for any $\e>0$ (see Corollary ~\ref{reg-lifts}).
\item If $g$ is $\e$-regular  on an  open  set $U\subset X$ then $g$ is $\e$-co-Lipschitz on $U$. (This is an easy corollary of Lemma~\ref{per-lem}. See~\cite{Per-Morse, Per-DC} for details.)
\item $g\co X_{reg}(g)\to \R^k$ is open. 
\item $k\le \dim X$ if  $X_{reg}(g)\ne \emptyset$).  This immediately follows from  part (1) of Lemma~\ref{per-lem}.
\end{enumerate}

\begin{thm}[Local Fibration Theorem]\label{fibration}\cite{Per-Morse, Per-DC}
Let $g\co X\to \R^k$ be an admissible map. Then $g|_{X_{reg}(g)}\co X_{reg}(g)\to \R^k$ is locally a topological bundle map. This means that any point $p\in  X_{reg}(g)$ possesses an open product neighborhood with respect to $g$ with an $MSC$-space as a fiber.
\end{thm}

We will use the original  Local Fibration Theorem without a proof. However, we will prove a more general version of it  (see Theorem  ~\ref{locrelfib}) in Section ~\ref{sec:extrstab} in order to prove the Relative Stability Theorem.

We will need the following  technical Lemma  due to Perelman which plays a key  role in the proofs of both the Stability Theorem and the Fibration Theorem above.

\begin{klem}\label{incompl}~\cite[Section 3]{Per-Morse}
 Let $p$ be a regular point of $f\co X\to \R^k$.  Suppose $f$ is incomplementable at $p$, i.e.  for any admissible function  $f_1\co X\to\R$, the point $p$ is not regular for $(f,f_1)\co X\to \R^{k+1}$. 

Then there exists an admissible function $h\co X\to\R$ with the following properties
\begin{enumerate}[(i)]
\item  $h(p)=0$.
\item $h$ is strictly concave on $B(p,R)$ for some $R>0$.
\item There are $r>0, A>0$ such that  $h< A$ on $f^{-1}\left(\bar{I}^k(f(p),r) \right)$ and   $f^{-1}\left(\bar{I}^k(f(p),r) \right)\cap  \{h\ge -A\}$ is compact in $B(p,R)$.

\item $h$ has a unique maximum in $B(p,R)\cap f^{-1}(v)$ for all $v\in \bar{I}^k(f(p),r) $.  Let $S$ denote the set of such maximum points.
\item $(f,h)$ is regular on $\left(\bar{I}^k(f(p),r) \cap B(p,R)\right)\backslash S$.
\end{enumerate}
\end{klem}

As an immediate corollary we obtain
\begin{cor}\label{bilip}
Let $f\co X^n\to\R^n$ be an admissible map. Then $f$ is locally bi-Lipschitz on the set of its regular points $X_{reg}(f)$.
\end{cor}

Let us first give an informal idea of the proof. Suppose we have $n+1$ points $a_0,\ldots, a_{n}$ such that  all the comparison angles at $p$ for  the triangles $\Delta a_ipa_j$ are  bigger than $\pi/2$. Then $f=(d(\cdot, a_1),\ldots, d(\cdot, a_{n}))$ is regular at $p$. Suppose $f$ is not injective near $p$. Then there exist points $x$, $y$ close to $p$ such that $f(x)=f(y)$. Let $z=\g(1/2)$ be the middle of a shortest geodesic $\g\co [0,1]\to X$ connecting  $x$ to $y$. Since the triangle $\Delta xa_iy$ is isosceles for every $i$, Toponogov comparison implies that
$\g'(1/2)$ is almost perpendicular to $\Uparrow_z^{a_i}$ for every $i=1,\ldots, n$ and the same is true for the vector opposite to $\g'(1/2)$.  On the other hand, we still have that $f$ is regular at $z$ and  $\angle \Uparrow_z^{a_i}\Uparrow_z^{a_j}>\pi/2$ for all $i\ne j$ provided $x, y$ are sufficiently close to $p$.  It's easy to see that this is impossible for dimension reasons which gives a contradiction.

\begin{proof}[Proof  of Corollary~\ref{bilip}]
Let $p$ be a regular point of $f$.  By part a) of Lemma~\ref{per-lem}, $f$ is incomplementable at $p$. Let $h$ be the function in $B(p,R)$ provided By lemma~\ref{incompl}.  By part (v) of  Lemma~\ref{incompl}, $(f,h)\co B(p,R)\backslash S\to\R^{n+1}$ is regular. However, as was just mentioned, by  part (1) of Lemma~\ref{per-lem},  a map from $X^n$ to $\R^{n+1}$ can not have any regular points. Therefore $B(p,R)\backslash S=\emptyset$. By part (iv) of Lemma \ref{incompl}, this is equivalent to saying that $B(p,R)\cap f^{-1}(v)$  consists of a single point for all $v\in \bar{I}^k(f(p),r)$  which means that $f$ is  injective near $p$. Finally, recall that being regular,  $f$ is both Lipschitz and co-Lipschitz near $p$, which together with local injectivity means that it's locally bi-Lipschitz on  $X_{reg}(f)$.
\end{proof}

\begin{proof}[Proof of Key Lemma~\ref{incompl}]
Since a complete proof is given in \cite{Per-Morse} we don't include as many details.

For simplicity we assume that all components  $f_j$ of $f$ are actually concave near $p$ (the proof can be easily adapted to the general case of semiconcave $f_j$s).  Since $f$ is $\e$-regular at $p$ there is a point $q$ near $p$ such that $f_j(q)>f_j(p)+\e d(p,q)$ for all $j=1,\ldots,k$. Since $f_j$'s are  $1$-Lipschitz, for all $r$ near $p$ and all $x\in B(q,\e d(p,q)/4)$ we have

\begin{equation}
f_j(x)>f_j(r)+\e \frac{d(x,r)}{4}
\end{equation}

Fix a small  positive $\d\ll \e d(p,q)$.  Choose a maximal $\d$-net $\{q_\a\}_{\a\in \mathfrak A}$ in $S(q,\e d(p,q)/4)$.  A standard volume comparison argument shows that $|\mathfrak A|=N\ge c\d^{1-n}$ where $n=\dim X$.

Consider the function $h=N^{-1}\Sigma_\a h_\a$ where $h_\a=\phi_\a(d(\cdot, q_a))$  and $\phi_\a\co \R\to \R$ is the unique continuous  function  satisfying the following properties

\begin{enumerate}
\item $\phi_\a'(t)=1$ for $t\le d(p,q_\a)-\d$;
\item $\phi_\a'(t)=1/2$ for $t\ge d(p,q_\a)+\d$;
\item $\phi_\a''(t)=-\frac{1}{4\d}$ for $d(p,q_\a)-\d<t<d(p,q_\a)+\d$;
\item $\phi_a(d(p,q_\a))=0$
\end{enumerate}

It is clear that $h$ is admissible near $p$ with $h(p)=0$ and

\begin{equation}\label{e:grad}
\langle d_rh,d_rf_j\rangle<-\e/8 \text{ for  all $j=1,\ldots, k$ and all $r$ near $p$}
\end{equation}

By ~\cite[Lemma 3.6]{Per-Morse}, $h$ is strictly $c\d^{-1}$-concave on $B(p,\d)$ ( see also~\cite[Lemma 4.2]{Kap} for a  more detailed proof of the same statement).
 
Denote $\Sigma_p^\e=\{\xi\in \Sigma_p \quad |\quad f_j'(\xi)>\e \text{ for all } j=1,\ldots,k\}$.

It easily follows from the definition that  $f$ is incomplementable at $p$ iff $\diam (\Sigma_p^0)\le \pi/2$.

\begin{slem}\label{estimate}
If $\diam (\Sigma_p^0)\le \pi/2$ then for all $r\in B_p(\d)$ we have

\begin{equation}\label{e:estimate}
h(r)\le h(p)-c\cdot d(p,r)+c\cdot max_j\{0,f_j(p)-f_j(r)\}
\end{equation}

\end{slem}

\begin{proof}
Since $h$ is $1$-Lipschitz and $f$  is $\e$-co-Lipschitz, by using Lemma~\ref{colip}, it's enough to prove the Sublemma for $r\in B(p,\d)$ satisfying $f_j(r)\ge f_j(p)$ for all $j=1,\ldots,k$. Then $\Uparrow_p^r\subset \Sigma_p^0$. By construction we also have that $\Uparrow_p^{q_\a}\subset \Sigma_p^0$ which by assumption of the sublemma implies that 
 that $\angle \Uparrow_p^r \Uparrow_p^{q_\a} \le \pi/2$ for all $\a$.  (This is the only place in the proof where we use that  $\diam (\Sigma_p^0)\le \pi/2$ i.e that $f$ is incomplementable at $p$). This means that the derivative of $h_\a$ at $p$ in the direction of $r$ is $\le 0$. By semi-concavity of $h_\a$ this implies
\begin{equation}\label{e:1}
h_\a(r)\le h_\a(p)-\la \cdot d(p,r)^2
\end{equation}

Moreover,  a volume comparison argument (see~\cite{Per-Morse} or ~\cite[Lemma 4.2]{Kap} for details)  shows that for most $\a\in \mathfrak A$ we actually have $\angle \Uparrow_p^r \Uparrow_p^{q_\a} \le \pi/2-C$ .

 Indeed, recall that $N\ge c\d^{1-n}$.

Fix a small $\mu>0$. By the first variation formula and semi-concavity of $h_\a$  we see that  if $\angle \Uparrow_p^r \Uparrow_p^{q_\a} \le \pi/2-\mu$ then  

\begin{equation}\label{b1}
h_\a(r)\le h_\a(p)-\la \cdot d(p,r)^2-\frac{1}{2}\mu \cdot d(p,r)
\end{equation}

A standard volume comparison argument shows that  the $n-1$-volume of the set 

$A_\mu=\{\xi\in\Sigma_p$ such that $ \pi/2-\mu \le \angle \xi  \Uparrow_p^r\le \pi/2\}$ is bounded above by $c\mu$.   By another standard volume comparison this implies that the maximal number of points in $A_\mu$ with pairwise angles  $\ge \d$ is at most $c\mu\d^{1-n}$.

This means that if $\mu\ll c$ then  for  the vast majority of $q_\a$ we must have $\angle \Uparrow_p^r \Uparrow_p^{q_\a} \le \pi/2-\mu$.

A suitable choice of $\mu$ now immediately yields $h(r)\le h(p)-c\cdot d(p,r)$.
\end{proof}

Sublemma~\ref{estimate} obviously implies that  inside $B(p,\d)$ the  sets $\{h\ge -c\d\}\cap f^{-1}(v)$ are compact for all $v$ close to $f(p)$ which proves (iii).

It  remains to prove parts (iv) and (v) of Key Lemma~\ref{incompl}.
For any $v\in \R^k$ denote $U_v=B(p,\d)\cap f^{-1}(v)$ and $U_v^+=\{x\in B(p,\d)\quad |\quad f_j(x)\ge v_j \text{ for all } j=1,\ldots,k\}$

\begin{slem}\label{max-fibers}
Let $z\in U_v$ be a point of maximum of $h$ on $U_v$ where $|f_j(p)-v_j|\le \d^2, j=1,\ldots, k$.
Then for every $x\in U_v\cap B(p,\d/2)$ we have
\[
h(z)\ge h(x)+c\d^{-1}d(x,z)^2
\]

\end{slem}

\begin{proof}

First we notice that
\begin{equation}\label{max}
\max_{U_v^+} h=\max_{U_v^+\cap \bar{B}(p,\d/2)} h=\max_{U_v \cap \bar{B}(p,\d/2)} h=h(z)
\end{equation}

The first and the last equalities follow from Sublemma~\ref{estimate} and the fact that $h$ is $1$-Lipschitz and $f$ is $\e$-co-Lipschitz.
If the equality in the middle is violated, then there exists a maximum point $r\in U_v^+\cap \bar{B}(p,\d/2)$ such that $f_l>v_l$ for some $l$. By applying part (2) of Lemma~\ref{per-lem} to $d_rh, d_r f_j$ in $\Sigma_rX$ we can find a direction $\xi\in \Sigma_rX$ such that $h'(\xi)>0$ and $f_j'(\xi)>0$ for all $j\ne l$. This contradicts the fact that $r$ is a point of maximum. This proves~(\ref{max}).

Now consider the midpoint $y$ of a shortest curve connecting $x$ and $z$. By concavity of $f$ and strict concavity of $h$ we get

\begin{equation}\label{p2}
h(y)\ge \frac{h(x)+h(z)}{2}+c\d^{-1}d(x,z)^2, \text{ and }
\end{equation}
\[
 f_j(y)\ge f_j(z) \text{ for all } j=1,\ldots,k
\]
In particular $y\in U_v^+$ and therefore $h(y)\le h(z)$. Combining this with (\ref{p2}) we immediately get the statement of Sublemma~\ref{max-fibers}.

\end{proof}
To complete the proof of Lemma~\ref{incompl} it remain to verify (v).

Let $z\in U_v\cap B(p,\d)$ be the point of maximum of $h$ on $U_v$. For any other point $x\in U_v$, by Sublemma~\ref{estimate}  combined with Lemma~\ref{colip}, we can find a point $s$ arbitrary close to $z$ and such that
\[
f_j(s)>f_j(x) \text{ for all } j \text { and } h(s)>h(x)+c\d^{-1}d(x,s)^2
\]

Let $\xi=\uparrow_x^s$.  By concavity of $h$ and $f_j$'s it is obvious that $h'(\xi)>0$ and $f_j'(\xi)>0$ for all $j$. Combined with (\ref{e:grad}) this means that $(f,h)$ is regular at $x$.

\end{proof}

\begin{rmk}
Key Lemma~\ref{incompl} in conjunction with Theorem~\ref{local bundle}  easily yields  Fibration Theorem~\ref{fibration} (see~\cite{Per-Morse} for details).

We will  prove a more general version of it in Section~\ref{sec:extrstab}.
\end{rmk}

We will need the following strengthened version of the Key Lemma which is the main geometric ingredient in the proof of the Stability Theorem.

\begin{lem}\label{incompl-lift}
Under the conditions of Lemma~\ref{incompl} suppose we have a noncollpasing converging sequence $X^n_i\convGH X^n$ and admissible functions $f_i$ on $X_i$ converging to $f$. Let $p_i\in X_i$  satisfy $p_i\conv p$.

Then there exist  admissible lifts $h_i$ of the function $h$ provided by Lemma~\ref{incompl} such that for all large $i$, the functions $f_i, h_i$ satisfy the properties (i)-(v) of Lemma~\ref{incompl}.
\end{lem}
\begin{proof}
Part (i) is obvious by construction of $h$ and $h_i$  where we might have to shift $h_i$ by a small constant to make it zero at $p_i$.

The fact that the natural  lifts $h_i$  of $h$ are strictly concave in $B(p_i,R)$  is basically the same as the proof of the concavity of $h$ itself (see  \cite[Lemma 3.6]{Per-Morse}). It is carried out in  full detail in~\cite[Lemma 4.2]{Kap}. This proves part  (ii).

For the proof of (iii) we can not use Sublemma~\ref{estimate}  as we did in the proof of Lemma~\ref{incompl} because $f_i$ might be complementable  everywhere near  $p_i$.  Nevertheless,   part (iii) is obvious because  it holds for $(f,h)$ by Lemma~\ref{incompl}  and   $f_i\to f, h_i\to h, p_i\to p$.

Parts (iv) and (v) are proved in exactly the same way as in the proof of Lemma~\ref{incompl}.

\end{proof}

\begin{lem}\label{reg-inner}
Let $F=(f_1,\ldots,f_k)\co X\to \R^k$ be $\e$-regular near $p$.   Then the level set $H=\{F={F(p)}\}$  is locally $\varkappa$-connected near $p$ for a linear function $\varkappa$.
\end{lem}

\begin{proof}
In what follows all constants $C$ will depend on $\e$.

Denote  $F(p)=v=(v_1,\ldots,v_k)$ and let
 $H_-$ be the set $\cap_i\{f_i\le v_i\}$.

By definition of a regular point, there exists $q$ near $p$ such that for $\xi=\uparrow_p^q$ we have $f_i'(\xi)>\e$ for all $i$ ( and the same holds for all $z$ near $p$).  Let $\eta=\nabla_p d(\cdot, q)$. Then $\eta$ is  polar  to $\xi$ by~(\ref{e:polarg}).  Therefore $df_i(\xi)+df_i(\eta)\le 0$ for all $i$  by ~(\ref{e:polar})  and hence

\begin{equation}\label{p3}
df_i(\eta)\le -\e \text{ for all } i.
\end{equation}

Let $x,y\in B(p,r)$ with $r\ll d(p,q)$ be two close points on the level set $H$. Let $\g$ be  a shortest geodesic connecting $x$ and $y$. Consider the gradient flow $\phi_t$ of $d(\cdot, q)$ in $B(p,r)$. 

By above, all $f_i$'s decrease with the speed at least $\e$ along $\phi_t$. Since all $f_i$ are 1-Lipschitz we know that along $\g$  we have
$f_i\le f_i(p)+ d(x,y)$.  By (\ref{p3}) this implies that for some $t\le C\cdot d(x,y)$ we can guarantee that  $\g_1=\phi_t(\g)$ lies in $H_-$.

By shifting $f_i$'s by constants we can assume that all $v_i$'s are the same and equal to $a$. Since $d(\cdot, q)$ is $\la$-concave 
in $B(p,r)$ ( where $\la$ depends on $d(p,q)$ and the lower curvature bound of $X$), 
by the Lipschitz properties of the gradient flows (\ref{e:lip}),  we know that $L(\g_1)\le C\cdot d(x,y)$.

Let $\g_2$ be the concatenation of the gradient curve of $d(\cdot, q)$ through $x$ followed by $\g_1$ followed by the  the gradient curve of $d(\cdot, q)$ through $y$  taken in opposite direction.  Then by above we still have
that $L(\g_2)\le C\cdot d(x,y)$.  We have also shown that $\g_2\subset H_-$.

Let $f=\min(a, \min_i f_i)$. Then $f$ is still semiconcave with the same concavity constants as $f_i$'s.

Observe that the gradient flow of $f$ takes $H_-$ to $H$ (the points of $H$ do not move under the flow).

By construction we have that $f\ge f(p)-C\cdot  d(x,y)$ along $\g_2$.
By $\e$-regularity of $F$ near $p$  and  Lemma~\ref{per-lem}   we see that $|\nabla f|\ge \e$  on  $\{f<a\}\cap B(p,r)$ and hence
 the gradient  flow of $f$ pushes $\g_2$ into a curve $\g_3 \subset H$  connecting $x$ and $y$ in  uniformly bounded time.  Once again applying the Lipschitz properties  of gradient flows we obtain  that $L(\g_3)\le C\cdot d(x,y)$.

\end{proof}

\begin{rmk}\label{cont-inner}
It easily follows from the proof  that the function $\varkappa$ provided by Lemma~\ref{reg-inner} is semi-continuous under Gromov-Hausdorff convergence in the following sense.

Suppose  $X_i^n\to X^n$ be a convergent sequence of {\it compact}  Alexandrov spaces with $\curv\ge \kappa$ and $f_i\co X_i\to \R^k$ be a sequence of admissible maps with 
 $\la$-concave  $1$-Lipschitz components  converging to $f\co X\to \R^k$. Suppose $f$ is $\e$-regular near $p$  and $p_i\in X_i$ converges to $p$. Then, by Corollary~\ref{reg-lifts},  $f_i$ is $\e$-regular on $B(p_i,r)$ for some $r>0$  for all large $i$ and the level sets $\{f_i=f_i(p_i)\}$ are $\varkappa$-connected near $p_i$ for all large $i$  with  the same $\varkappa(t)=C(\e,\kappa,\la, X)\cdot t$.
\end{rmk}

\section{Proof of the stability theorem}\label{sec:proofstability}
\begin{defn}
Let $p$ be a point in an Alexandrov space $X$.  Let $g=G\circ f=(g^1,\ldots,g^k)\co X\to \R^k$ be regular  at $p$.
Then  an open  product neighborhood of $p$  with respect to such $g$ is called 
a product neighborhood of a $p$ of rank $k$.
\end{defn}
\begin{defn}
We'll call a subset $H$ of $\R^n$ a {\sf generalized quadrant} if it has the following form $$H=\cap_{i\in I}\{x_i\ge c_i\}\cap_{j\in J}\{x_j\le c_j\}$$
where $I,J$ are some (possibly empty) subsets of $\{1,\ldots,n\}$ and $(x_1,\ldots,x_n)$ are the standard coordinates on $\R^n$.

\end{defn}
\begin{defn}
A compact subset $P$ in an Alexandrov space $X$ is called $k$-framed if 
 $P$ can be covered by a finite collection of open sets $U_\a$
 such that each $U_\a$ is a product neighborhood of rank $k_\a \ge k$ for some $p_\a\in P$  with respect to some  $$g_\a=G_\a\circ f_\a\co X\to \R^{k_\a}$$ and $P\cap U_\a=g_\a^{-1}(H_a)\cap U_\a$  where $H_\a$ is a generalized quadrant in $R^{k_\a}$. 
 
 We say that a framing $\{U_\a, f_\a, H_\a\}_{\a\in\mathfrak A}$ respects a map $f\co X\to\R^l$ if the first $l$ functions of $f_\a$ coincide with $f$ for all $\a$.

 More generally, we'll  say that a framing $\{U_\a, g_\a=G_\a\circ f_\a, H_\a\}_{\a\in\mathfrak A}$ respects a representation of an admissible map $g$  $$g=G\circ f\co X\to \R^l$$ if 
 the first $l$ functions of each $f_\a$ coincide with $f$ and $G_\a$ has the form $G_\a(x,y)=(G(x), T_\a(x,y))$  where $x\in \R^l, y\in \R^{k_\a-l}$, i.e the first $l$ coordinates of $G_\a$ are equal to $G$. (In particular, the first $l$ functions of $g_\a$ also coincide with $g$).

 \end{defn}

To simplify notations,  we'll often simply say that a framing respects a map $g$ to mean that it respects the representation $g=G\circ f$.

\begin{ex}
Any compact Alexandrov space has a zero framing.
\end{ex}

\begin{ex}
For any point $p\in X$ there exists $\e>0$ such that $d(\cdot, p)$ has no critical points in $B(p,\e)\backslash \{p\}$. Hence, for any positive $r<R<\e$, the annulus $A(r,R,p)=\bar{B}(p,R)\backslash B(p.r)$ is 1-framed respecting $f=d(\cdot, p)$.
\end{ex}

Suppose $X^n_i\to X^n$ is a converging noncollapsing sequence of compact Alexandrov spaces with $\curv\ge \kappa$. Let
$\theta_i\co X\to X_i$ be a sequence of $o(i)$-Hausdorff approximations.

Let $P\subset X$ be a $k$-framed compact subset of $X$. We define the corresponding $k$-framed subsets $P_i\subset X_i$ as follows.   Let $g_\a=G_\a\circ f_\a$ be a representation of $g_\a$ given by the definition of an admissible map. 
We  lift  the defining functions $f_\a$ and $f$ to $f_i$ and $f_{\a,i}$ in the natural way.  Suppose $G_{\a,i}$ is  a sequence of {\it uniformly}  bi-Lipschitz homeomorphisms of open sets in $\R^k$ converging to $G$. Put $g_{\a,i}=G_{\a,i}\circ f_{\a,i}$.

Then $g_{\a,i}$ will still be admissible and regular on the corresponding subsets of $X_i$ by Corollary~\ref{reg-lifts}.  In particular we get product neighborhoods $U_{\a,i}$ with respect to $g_{\a,i}$.  


We'll say that a compact set $P_i\subset X_i$ is a lifting of $P$ if  $P_i\cap U_{\a,i}=g_{\a,i}^{-1}(H_a)$ for all $\a$.

\begin{rmk}
Note that a lifting of $P$ need not exists! However, if it does, it is automatically $k$-framed. Moreover, if $P$ is $k$-framed with respect to $f\co X\to \R^l$ and $I_\a,J_\a\subseteq\{1,\ldots, l\}$ for all $\a$ then the lifting exists for all large $i$.
In particular, if $X$ is a compact Alexandrov space and $X^n_i\to X^n$ with $\curv\ge \kappa$ then the lifting of $X$ with respect to a zero framing exists for any large $i$ and is equal to $X_i$.
\end{rmk}

\begin{lem}\label{k-frammed-inner}
The set's $P_i$ are $\varkappa$-connected for all large $i$ and the same $\varkappa(t)=C t$.
\end{lem}
\begin{proof}
Let $x,y$ be two close points in $P\cap U_\a$ for some $\a$. Since $g_\a$ is L-Lipschitz, $|g_\a(x)-g_a(y)|\le L\cdot d(x,y)$. Since $H_a$ is convex, the straight line segment connecting $g_\a(x)$ and $g_\a(y)$ lies in $H_\a$.  Since  $f_\a$ (and hence$g_\a$) is $\e$-co-Lipschitz , by Lemma~\ref{colip} we can lift it to a curve $\g_1\co [0,1]\to U_\a\cap P$ of length $\le \frac{L}{\e}d(x,y)$ with $\g_1(0)=x$. 
Observe that $\g_1(1)$ and $y$ lie in the same fiber of $g_\a$ (and hence of $f_\a$) and $d(\g_1(1),y)\le C(L,\e)d(x,y)$ by the triangle inequality.
By Lemma~\ref{reg-inner} we can connect $\g_1(1)$ and $y$ by a curve $\g_2$ inside the  fiber $\{g_\a=g_\a (y)\}$ of length $\le \tilde{C}\cdot d(\g_1(1),y)$. The concatenation of $\g_1$ and $\g_2$ provides a curve $\g$ in $P\cap U_\a$ connecting $x$ to $y$ with $L(\g)\le  C\cdot d(x,y)$.

As was observed in Remark~\ref{cont-inner} the constant $\tilde{C}$ in the above argument can be chosen to be the same for all $f_{\a,i}$ and hence, since all $G_{\a,i}$ are {\it uniformly}  bi-Lipschitz, all $P_i$ are $\varkappa$-connected for all large $i$ for the same $\varkappa(t)=C t$.

\end{proof}

The proof of the stability theorem proceeds by reverse induction in framing and, in fact, it requires us to
 to prove  the following stronger version of it:

\begin{thm}[Parameterized Stability Theorem]\label{gen-stab}
Suppose $X^n_i\to X^n$ is  a converging noncollapsing sequence of Alexandrov spaces with curvature bounded below and diameter bounded above. Let
$\theta_i\co X\to X_i$ be a sequence $o(i)$-Hausdorff approximations.

Let $P\subset X$ be a $k$-framed compact subset of $X$ whose framing respects  $g\co X\to\R^l$.  Let $K\subset P$ be a compact subset such that the framing of $P$ respects $g$ on $P$ and $(g,h)$ on $K$ for some $g\co X\to \R^{l}$ and $h\co X\to \R$.

Then for all large $i$ there exist homeomorphisms $\theta_i'\co P\to P_i$ such that
$\theta_i$ is $o(i)$-close to $\theta_i$ and  respects $g$ on $P$ and $(g,h)$ on $K$.
\end{thm}

\begin{proof}
We proceed by reverse induction in $k$.

If $k=n$ then the locally defined maps $g_{\a,i}\co U_{\a,i}\to \R^n,  g_\a\co U_\a\to\R^n$ are bi-Lipschitz homeomorphisms. By construction, the maps $\theta_{\a,i}=g_{\a,i}^{-1}\circ g_\a\co U_\a\to U_{\a,i}$ are homeomorphisms Hausdorff  close to $\theta_i$. Moreover, by construction, 
$\theta_{\a,i}$ sends $P\cap U_\a$ onto $P_i\cap  U_{\a,i}$.

Thus the statement of the theorem follows from Strong Gluing Theorem~\ref{stronggluing} and Lemma~\ref{k-frammed-inner}.

Induction step.

Suppose the theorem is proved for $k+1\le n$ and we need to prove it for $k$.

Let $P$ be $k$-framed and let $p$ lie in $P$. Then $p\in U_\a$ for some $\a$.  Let $g_\a=G_\a\circ f_\a\co U_\a\to \R^k$ be the admissible map regular at $p$ coming from the definition of a $k$-framed set.  

To simplify the notations we will assume that $G_\a=\Id$ and $f_\a=g_\a$. The proof in the general case easily follows from this one with obvious modifications.

 Let $p_i=\theta_i(p)$.

By possibly adding more components to $f_\a$ we can construct an  admissible  map to $f_p\co X\to\R^{k_p}$ where $k_p\ge k$ which is incomplementable at $p$. 

 Let $h\co B(p,R)\to R$ be a strictly concave function provided by Lemma~\ref{incompl}. By choosing a sufficiently small $r, A>0$ we can assume that the set
$U_p=f_p^{-1}(I^{k_p}(f_p(p),r)\cap \{h\ge -A\}$ is compact. By reducing $r$ further we can assume that $|h(x)|\le a\ll A$ for $x\in S\cap U_p$. We will call $U_p$  a special neighborhood of $p$.

By Lemma~\ref{incompl} we have that $U_p$ is a $k$-framed compact subset of $X$.

Since $h$ and all the  coordinates of $f_p$  are  admissible, they have natural  admissible lifts $h_i$ and $f_{p,i}$ which define  corresponding neighborhoods $U_{p_i}$ of $p_i$.

The proof of  Stability Theorem~\ref{gen-stab}  will easily follow from the following 

\begin{lslem}\label{keylemma}
For all large $i$ there exist homeomorphisms $\theta_{p,i}\co U_p\to U_{p_i}$ respecting $f_p$ and $o(i)$-close to the Hausdorff approximation $\theta_i$.
\end{lslem}

Let us first explain how to finish the proof of theorem~\ref{gen-stab}  given Lemma~\ref{keylemma}.

Choose a finite cover of $P$ by the interiors of the special neighborhoods $U_{p_\beta}$. For all large $i$,  Lemma~\ref{keylemma} provides homeomorphisms $\theta_{p_\beta,i}\co  U_{p_\beta}\to  U_{p_\beta,i}$ respecting $f_{p_\beta}$   and $o(i)$-close to $\theta_i$.

Observe that since each $\theta_{p_\beta,i}$ respects $f_{p_\beta}$, it sends $P\cap  U_{p_\beta}$  onto $P_i\cap U_{p_\beta,i}$.

Taking into account Lemma~\ref{k-frammed-inner} we can  apply  Gluing Theorem~\ref{stronggluing} to obtain the desired homeomorphism $\theta_i'\co P\to P_i$.

\begin{proof}[Proof of Local Stability Lemma~\ref{keylemma}]

All throughout the proof of the Lemma we will work only with points in $U_p$  in $X$ and $U_{p_i}$ in $X_i$.

If $k_p>k$ then the statement follows directly from the induction hypothesis.

Let's suppose $k_p=k$.

First we change the function $h$ to an auxiliary function $\th$ by shifting $h$ by a constant on each fiber of $f$ to make it identically zero on $S$.

More precisely, let $\th(x)=h(x)-h(S\cap f^{-1}(f(x)))$. Recall that  by  Lemma~\ref{incompl}(iii),  $S\cap f^{-1}(f(x)))$ consists of a single point so that this definition makes sense.  Also by Lemma~\ref{incompl}(iii),   we have $h(S\cap f^{-1}(f(x)))=\max_{y\in U_p\cap  f^{-1}(f(x))} h(y)$ and therefore
 $\th=h-H\circ f$ where $H\co \R^k\to\R$ is  given by $H(v)=\max_{x\in U_p\cap f^{-1}(v)}h(x)$. Since $f$ is co-Lipschitz and $h$ is Lipschitz , using Lemma~\ref{colip} we easily conclude  that $H$ is Lipschitz.
 

 In particular, $(f,\th)=\bar{H}\circ (f, h)$ where $\bar{H}$ is a bi-Lipshitz homeomorphism of some open domains in $\R^{k+1}$ given by $\bar{H}(a,b)=(a, b-H(a))$. 

Therefore, we still have that $(f, \th)$ is regular on $U_p\backslash S$ and hence it's locally a bundle map on $U_p\backslash S$ by Theorem~\ref{fibration}.

In addition, by construction, $\th=0$ on $S$ and $\th<0$ on $U_p\backslash S$.

We define $\th_i, H_i$ and $\bar{H}_i$  in a similar fashion using $f_i, h_i$.  We obviously have that $\th_i\to \th, H_i\to H, \bar{H}_i\to\bar{H}$. Moreover, since $f_i$ are uniformly co-lipschitz,  all $H_i$ are uniformly Lipschitz and hence  all $ \bar{H}_i$ are uniformly bi-Lipschitz.

Then we again have that $\th_i=0$ on $S_i$ and $\th<0$ on $U_{p_i}\backslash S_i$.  By Lemma~\ref{incompl-lift} we also have that
$(f_i,\th_i)$ is regular on $U_{p_i}\backslash S_i$.

Fix a small  $\d\ll A$. Then the set $\{\th \le -\d\}$ is $(k+1)$-framed  with the corresponding sets in $U_{p_i}$ given by $\{\th_i \le -\d\}$  and therefore, by induction assumption, for large $i$ there exist homeomorphisms $\theta_{\d,i}'\co \{\th \le -\d\}\to \{\th_i \le -\d\}$ respecting $(f,\th)$ and $o(i|\d)$-close to $\theta_i$.

In particular, the fiber $F_i$ of $(f_i, \th_i)$ is homeomorphic to the fiber $F$ of $(f,\th)$ for all large $i$.

Next  consider the set $\{\th> -3\d/2\}$ and consider the map $(f,\th)\co \{\th> -3\d/2\}\backslash S\to I^k(f(p),r)\times (-3\d/2,0)$.

By Lemma~\ref{incompl}, this map is regular and proper. Therefore, by Theorem~\ref{fibration} and Theorem~\ref{local bundle},
it is a fiber bundle.
 
Hence,  $\{\th> -3\d/2\}\backslash S$ is homeomorphic to $I^k(f(p),r)\times (-3\d/2,0)\times F$  with the first two coordinates given by $(f,\th)$.
By restriction this gives a homeomorphism  $\{\th\ge \d\}\backslash S$  to $I^k(f(p),r)\times [-\d,0)\times F$  with the first two coordinates still  given by $(f,\th)$.

By Lemma~\ref{incompl} and construction of $\th$, $\th$ has a unique max ( equal to zero) on $f^{-1}(v)$ for any $v\in  \bar{I}^k(f(p),r)$. Therefore, the above homeomorphism can be uniquely extended to a homeomorphism
$\phi_\d\co  \{\th\ge -\d\}\to \bar{I}^k(f(p),r)\times KF$ respecting $(f,\th)$. 

Similarly, for all large $i$ we have homeomorphisms $\phi_{\d,i}\co \{\th_i\ge -\d\}\to \bar{I}^k(f(p),r)\times KF_i$ respecting $(f_i,\th_i)$.

Recall that $F_i$ is homeomorphic to $F$ for large $i$. Modulo that homeomorphism we can take
 the composition $\phi_{\d,i}^{-1}\circ\phi_\d$ and  obtain a homeomorphism  $\theta_{\d,i}'' \co  \{\th\ge -\d\}\to  \{\th_i\ge -\d\}$ respecting $(f,\th)$.

Since $\theta_{\d,i}'' $ respects $f$ and $S_i$ is $o(i)$-close to $S$,  it is clear that $\theta_{\d,i}'' $ is $\varkappa(\d)+o(i|\d)$-close to $\theta_i$.

Gluing $\theta_{\d,i}'' $ and $\theta_{\d,i}' $ we obtain homeomorphisms  $\theta_{\d,i} \co U_p\to U_{p_i}$ respecting $f$ and $\varkappa(\d)+o(i|\d)$-close to $\theta_i$.

Since $\d$ was arbitrary,  a standard diagonal argument
provides  the desired homeomorphisms $\theta_{p,i}'\co U_p\to U_{p_i}$ respecting $f$ and $o(i)$- close to $\theta_i$.

\end{proof}
To conclude the proof of Theorem~\ref{gen-stab} observe that all the maps $f_p$  above can be chosen to respect $f$ on $P$ and $h$ on $K$.
\end{proof}

\begin{rmk}
It's instructive to point out precisely what's needed to make the proof of the Stability Theorem work in the Lipschitz category.

{\bf 1.}  One needs to generalize the deformation of homeomorphisms results of Siebenmannform~\cite{Sieb} used in the proof of Theorem~\ref{gluingthm}
to Lipschitz category. This is probably possible and in fact it is already known in case of Lipschitz manifolds by~\cite{Sul2}.

{\bf 2.} Another ( probably quite difficult ) point is to generalize Perelman's Local Fibration Theorem~\ref{fibration} to Lipschitz category. To do this one needs to show that under the assumptions of Lemma~\ref{incompl}, the homeomorphism of the "tubular" neighborhood of $S$   to the product of $S$ and the cone over $F$ can be made to be bi-Lipschitz.  The basic (and probably the most important) case of this would be to show that if $h$ is a proper strictly concave function on an Alexandrov space with a unique maximum at a point $p$ then the superlevel set $\{h\ge h(p)-\e\}$ is bi-Lipschitz homeomorphic to the cone over $\{h=h(p)-\e\}$. This is related to {\bf 1.}   and could possibly be proved using an appropriate generalization of Siebenmann's results.

Similar discussion applies to the case of DC rather than Lipschitz stability.

\end{rmk}

The  Stability theorem has a natural generalization to the case of  pointed Gromov-Hausdorff convergence.
The following application, saying that the stability homeomorhisms can be constructed on arbitrary large compact subsets,
 seems to be the most useful. For simplicity we only state the unparameterized version.

\begin{thm}\label{pointedstab}
Let $(X_i^n, p_i)\convGH (X^n,p)$ where $\curv X_i\ge \kappa$ for all $i$. Let $R>0,\e>0$  and let $\phi_i\co B(p,R+\e)\to B(p_i,R+\e)$ be $o(i|R)$-Hausdorff approximations. Then for all large $i$ there exist open embedding $\psi_i\co B(p,R+\e/2)\to X_i$ which are $o(i|R,\e)$-close to $\phi_i$ and such that $\psi_i(B(p,R+\e/2))\supset B(p_i,R)$ for all large $i$.
\end{thm}
\begin{proof}
Put $U=B(p,R), V=B(p,R+\e/2), W=B(p,R+\e)$. Recall, that for any point $x\in X$ there is $r_x>0$ such that $d(\cdot, x)$ has no critical points in $B(p,2r_x)\backslash \{p\}$. Cover $\bar{W}$ by finitely many such balls  $B(p_\a,r_\a)$ with all $r_\a<\e/4$.
Let $p_{\a,i}=\phi_i(p_\a)$. By  Stability Theorem~\ref{gen-stab}, for all large $i$  and all $\a$ there exist homeomorphisms $\phi_{\a,i}\co B(p_\a,r_\a)\to B(p_{\a,i},r_\a)$ which are $o(i|R,\e)$-close to $\phi_i$. Now the statement of the theorem follows by the direct application of Theorem~\ref{locgluing}.

\end{proof}
\begin{rmk}
In the proof of Theorem~\ref{pointedstab} we could not apply the stability theorem directly to $\bar{B}(p,R)$ because in general we have no information on the existence of critical points of $d(\cdot, p)$ outside a small ball around $p$. In particular, $\bar{B}(p,R)$ need not be an $MCS$-space or be $\varkappa$-connected for any large $R$.
\end{rmk}

\begin{rmk}
It follows from the stability theorem that for any given $n\in Z_+, k\in \R, D>0, v>0$ there exists $\e>0$ such that any $X$ and $Y$ in  $\Alkdv$ with $d_{G-H}(X,Y)\le \e$ are homeomorphic.
It's interesting to see if one can give an explicit estimate on $\e$ it terms of $n,k,D,v$.
\end{rmk}

\section{Finiteness of submersions}\label{sec:subm}
The following generalization of the Grove-Petersen-Wu finiteness theorem was proved by K. Tapp~\cite[Theorem  2]{Tapp}:
\begin{thm}\label{tappthm}
Given $n,k\in \Bbb Z_+, v,D,\lambda\in \Bbb R$ with $k\ge 4$, there are at most finitely many topologically equivalence classes of bundles in the set of Riemannian submersions $M^{n+k}\to B^n$ satisfying ${\rm vol}(B)\ge v, |{\rm sec}(B)|\le \lambda$; ${\rm vol}(M)\ge v, {\rm diam}(M)\le D, {\rm sec}(M)\ge -\lambda$.
\end{thm}

Here two submersions $\pi_i\co M_i\to B_i$ ($i=1,2$) are called topologically equivalent if there exist homeomorphisms $\phi\co M_1\to M_2, f\co B_1\to B_2$ such that $\pi_2\circ \phi=f\circ \pi_1$.

This definition can be naturally extended to the class of submetries of Alexandrov spaces.  
Recall that a map $f\co X\to Y$ is called a submetry if $f(B(x,r))=B(f(x),r)$ for any $x\in X, r>0$ (i.e if $f$ is both 1-Lipschitz and 1-co-Lipschitz).

It is obvious  that a Riemannian submersion between complete
Riemannian manifolds is a submetry. Moreover, converse is also
true according to \cite{BG}. It is also clear that if a compact
group $G$ acts on a Riemannian manifold $M$ by isometries, then
the projection $M\to M/G$ is a submetry.

Submetries enjoy many properties of Riemannian submersions. In particular, one can talk about horizontal and vertical tangent vectors and curves. Also, it's easy to see~\cite{BGP} that  submetries increase Alexandrov curvature, that is, if $\curv X\ge \kappa$ and $f\co X\to Y$ is a submetry, then $curv Y\ge \kappa$. For more basic information on submetries see~\cite{Lyt}.

Suppose $\pi\co X\to B$ is a submetry between compact Alexandrov spaces. 
It trivially follows from the definition that if $f\co B\to\R^k$ is   admissible then $f\circ\pi$ is  admissible on $M$. Moreover,
$f$ is regular at $p\in B$ iff $f\circ\pi$ is regular at any $y\in \pi^{-1}(p)$. In particular, if $B$ is a Riemannian manifold (or more generally, if $B$ is everywhere $n$-strained)  then $\pi$ is a fiber bundle. Thus the  above notion of equivalence of submersions naturally extends to submetries.

 Theorem~\ref{tappthm} generalizes a theorem of J. Y. Wu,~\cite{wu-fin} which proved the same result under a strong extra assumption that fibers of the submersions are totally geodesic. The proof of Theorem~\ref{tappthm}  relies on the proof of Wu's theorem  which just as the proof of  Grove-Petersen-Wu finiteness Theorem~\ref{GPW}  uses techniques of controlled homotopy theory. This explains the assumption $k\ge 4$ in Theorem~\ref{tappthm}. However, this assumption is, in fact, unnecessary as this result follows from the Parameterized Stability Theorem which does not require any dimensional restrictions.

\begin{thm}\label{sub-fin}
Given 
 $k\in \Bbb Z_+, V,D,\in \Bbb R_+, \la, \kappa\in \R$ , there are at most finitely many  equivalence classes   of   submetries $X^{n+k}\to B^n$  where $X^{n+k}\in  \text{{ \nnn Alex}}^{n+k}({D, \kappa,v})$  and $B^n$ is a closed Riemannian manifold satisfying ${\rm vol}(B)\ge V, |{\rm sec}(B)|\le \lambda$.
\end{thm}

\begin{proof}

We first give a proof in case of fixed $B$.

Let $\pi\co X\to B$ be a submetry where $X$ is an Alexandrov space of curvature bounded below.

Let $p\in B$ be any point.  Choose $n+1$ unit vectors $v_0,\ldots,v_n\in T_pB$ with pairwise angles bigger than $\pi/2$. Then for all sufficiently small $R>0$ the points $p_i=exp_p(Rv_i)$ define an admissible map $f\co B\to \R^n$ given by $y\mapsto (d(y,p_1),\ldots,d(y,p_n))$. This map is obviously regular on $B(p,r)$   for $r\ll R$  and  it gives a bi-Lipschitz  open embedding ( in fact a smooth one)  $f\co B(p,R)\to \R^n$.
Let $F_i=\pi^{-1}(p_i)$. Let $\tilde{f}\co X\to\R^n$ be given by $\tilde{f}(x)=(d(x,F_1),\ldots, d(x,F_n))$.

 Since $\pi$ is a submetry we obviously have that $d(x, F_i)=d(\pi(x),p_i)$ for any $i$ and any $x\in X$. Therefore $\tilde{f}\equiv f\circ\pi$.  It is also obvious that $\tilde{f}$ is regular at $x\in X$ iff $f$ is regular at $\pi(x)$. In particular, $\tilde{f}$ is regular on the $r$-neighborhood of  $\pi^{-1}(\pi(p))$. Thus, up to a bi-Lipschitz change of coordinates on the target, when restricted to $U_r(\pi^{-1}(\pi(p)))$, we can write $\pi$ as a proper regular map $\tilde{f}$ to $\R^n$. Let's cover $B$ by finitely many coordinate neighborhoods $U_\a=B(p_\a,r_\a)$  as above and let $f_\a\co U_\a\to\R^n$ be the corresponding coordinate projections. Since all $f_\a$ are bijections we obviously have that for any $x,y\in \pi^{-1}(U_\a\cap U_b)$, $\tilde{f}_\a(x)=\tilde{f}_\a(y)$ iff $\tilde{f}_\b(x)=\tilde{f}_\b(y)$.
 
 Therefore Parameterized Stability Theorem~\ref{gen-stab} can be easily amended to include the case when a framing on $X$ respects a submetry to a fixed manifold as all the arguments can be made local on $B$ where instead of the submetry $\pi$ one can work with  regular maps  $\tilde{f}_\a$.
 
 The case of variable base easily follows given the fact that by Cheeger-Gromov compactness  the class of manifolds  $\{B^n\quad |\quad 
 {\rm vol}(B)\ge V, |{\rm sec}(B)|\le \lambda, \diam(B)\le D\}$ is precompact in Lipschitz topology and its limit points are $C^{1,\a}$-Riemannian manifolds~(see e.g. ~\cite{GrWu}).
 
\end{proof}
\begin{rmk}
It's interesting to see whether Theorem~\ref{sub-fin} remains true if one removes the assumption about the uniform upper bound on the curvature of $B$.
\end{rmk}

\section {Stability with extremal subsets}\label{sec:extrstab}
The results proved in  this section are new.

The notion of an extremal subset in an Alexandrov space was introduced in ~\cite{Per-Pet}.
\begin{defn}
A closed subset $E$ in an Alexandrov space $X$ is called extremal if for any $q\in X\backslash E$ and $f=d(\cdot, E)$ the following holds:

If $p\in E$ is a point of local minimum of $f|_E$ then it's a critical point of maximum type of $f$ on $X$ , i.e.
\[
df_p(\xi)\le 0 \text{ for any } \xi \in \Sigma_pX
\]
\end{defn}

Alternatively, it was shown in ~\cite{Petrunin-survey}  that $E$ is extremal iff it's invariant under gradient flows of all semiconcave functions on $X$.

An extremal subset is called {\it primitive} if it doesn't contain any proper extremal subsets with nonempty relative interiors.

We refer to~\cite{Per-Pet,Petrunin-survey} for  basic properties of extremal subsets.
It is easy to see~\cite{Per-Pet}  that closures of topological strata in an Alexandrov space $X$  are extremal. Therefore stratification into extremal subsets can be considered as a geometric refinement of the topological stratification of $X$.

 It is of course obvious that a homeomorphism between two Alexandrov spaces has to preserve topological strata.

The goal of this section is to generalize the Stability Theorem by showing that the stability homeomorphisms can be chosen to preserve extremal subsets.
Namely we will prove the following

\begin{thm}[Relative Stability Thorem]\label{stab-extr}
Let $X^n_i\convGH X^n$ be a noncollapsing sequence of compact Alexandrov spaces in $\Alkd$.
Let $\theta_i\co X\to X_i$ be a sequence $o(i)$-Hausdorff approximations.
Let $E_i\subset X_i$ be a sequence of  extremal subsets converging to an  extremal subset $E$ in $X$.
Then for all large $i$ there exist homeomorphisms $\theta_i'\co (X,E)\to (X_i,E_i)$, $o(i)$-close to $\theta_i$.
\end{thm}

In order to prove this theorem we'll need to generalize all the machinery used in the proof of the regular Stability Theorem to its relative
version  respecting extremal subsets. This is fairly straightforward and only minor modifications of the proofs are required.

In particular we'll have to prove the relative version of Local Fibration Theorem~\ref{fibration}. Along the way we'll obtain some new topological information about the way a general extremal subset is embedded into an ambient Alexandrov space.

It was shown in ~\cite{Per-Pet} that just as Alexandrov spaces, extremal subsets are  naturally stratified  in the sense of the following definition.

\begin{defn}\label{MCS-t}
A metrizable space $X$ is called an $\widetilde{MSC}$-space of dimension $\le n$ if  every point $x\in X$ has a neighborhood pointed homeomorphic to an open cone over a compact   $\widetilde{MCS}$-space of dimension $\le n-1$. As for MCS-spaces  we assume the empty set to be the unique  $\widetilde{MCS}$-space of $\dim\le -1$.

\end{defn}

We will also need a relative version of the above definition.

\begin{defn}
Let $X$ be  an  $\widetilde{MSC}$-space of $\dim\le n$. A subset $E\subset X$ is called a stratified subspace of $X$  of dimension $\le k$ if every $x\in X$ has  a pointed neighborhood $(U,p)$ such that $(U,U\cap E,p)$ is homeomorphic to $(K\Sigma, K \Sigma', o)$ where $\Sigma$ is a compact $\widetilde{MSC}$-space of $\dim\le n-1$ and $\Sigma'\subset\Sigma$ is a compact stratified subspace of  $\Sigma$ of dimension $\le k-1$. As usual, the only stratified subspace of $\dim\le -1$ is the empty set.
\end{defn}

It is obvious from the definition that a stratified subspace in $X$ of $\dim\le k$ is an $\widetilde{MCS}$-space of dimension $\le k$. 

\begin{rmk}
It is easy to see that a connected  $\widetilde{MSC}$-space  is an $MCS$-space iff  its local topological dimension is constant.
It was shown in~\cite{Per-Pet} that  a primitive  extremal subset is equal to the closure of its top dimensional strata and therefore is an $MCS$-space by above.
\end{rmk}

In the process of proving the Relative Stability Theorem we'll obtain the following result which clarifies the relative topology of extremal subsets with respect to their ambient spaces.
\begin{thm}[Relative Stratification Theorem]\label{str-subspace}
Let $X$ be an Alexandrov space and let $E\subset X$ be an extremal subset. Then $E$ is a stratified subspace of $X$.
\end{thm}

Just as in the non-relative case, this theorem is a Corollary of the following local fibration theorem applied to the natural  map  $X\to \R^0$.

\begin{thm}[Relative Local Fibration Theorem]\label{locrelfib}
Let $f\co X\to \R^k$ be regular at $p\in E$ where $E\subset X$ is an extremal subset. Then there exists an open neighborhood $U$ of $p$, an $MCS$-space $A$, a stratified subspace $B\subset A$  and a homeomorphism $\phi\co  (U,E\cap U)\to (A,B)\times \R^k$ such that $\pi_2\circ \phi=f$.
\end{thm}

It was shown in~\cite{Per-Pet} by Perelman and Petrunin  that the intrinsic metric on an extremal subset of an Alexandrov space is locally bi-Lipschitz to the ambient metric. On closer examination their proof actually gives the following somewhat stronger statement which we'll need for the proof of the relative stability theorem:
\begin{lem}\label{uniform-bi-Lip}
There exists $\e=\e(n. D, \kappa, v)>0$ such that if $X\in \Alkdv$ and $E\subset X$ is extremal, then for any $p,q\in E$ with $d(p,q)\le \e$ there exists a curve in $E$ connecting $p$ and $q$ of length $\le \e^{-1}d(p,q)$.
\end{lem}
\begin{proof}
Because the argument is very easy we give it here for reader's convenience.
It is well-known (see ~\cite{Per-Pet} or~\cite{Gr-Pet}) that for the class $\Alkdv$ there exists an $\e>0$ such that the following holds:

If $X\in  \Alkdv$ and $p,q\in X$ with $d(p,q)< \e^2$ then
\begin{equation}\label{doublecrit}
|\nabla_p d(\cdot, q)|>\e \text{ or } |\nabla_q d(\cdot, p)|>\e
\end{equation}

Suppose the first alternative holds. Then there exists $x$ near $p$ such that $\tangle xpq\ge \pi/2+\e$.   Then $\nabla_p d(\cdot, x)$ is polar to $\Uparrow_p^x$ so that $\angle \nabla_p d(\cdot, x)\uparrow_p^q\le \pi/2-\e$.  This means that moving $p$  along the gradient flow of $d(\cdot, x)$ decreases $d(p,q)$ in the first order (with the derivative at $0$  at least $\e$). Since $E$ is extremal, the flow through $p$ remains in $E$.

 Now a standard argument shows that we can construct a curve in $E$  connecting $p$ and $q$ of length $\le \e^{-1}d(p,q)$.

\end{proof}

We will need the following generalization of Lemma~\ref{per-lem}
proved in~\cite{Per-Pet}:

\begin{lem}\cite{Per-Pet,Petrunin-survey}\label{per-lem-extr}
Suppose $\Sigma^{n-1}$ has $\curv\ge 1$ and let $f_0,\ldots, f_k\co \Sigma\to\R$ be functions of class $DER$  such that
$\e=\min_{i\ne j}-\langle f_i,f_j\rangle  >0$. Let $E\subset \Sigma$ be an extremal subset.
Then

\begin{enumerate}
\item There exists $w\in E $ such that $f_i(w)> \e$ for all $i\ne 0$.
\item There exists $v\in E$ such that $f_0(v)>\e, f_1(v)<-\e$ and $f_i(v)=0$ for $i=2,\ldots,k$.
\end{enumerate}
\end{lem}

Just as in the  case of regular functions on Alexandrov spaces, this lemma implies that if $f\co X\to\R^k$ is regular at $p\in E$ where $E\subset X$ is extremal then $f|_E$ is co-lipschitz near $p$. 

The main geometric ingredient in the proof of the Relative Local Fibration Theorem  and the Relative Stratification Theorem is the following relative analogue of Lemma~\ref{incompl}.

\begin{lem}\label{relincompl}
 Let $p\in E$ be a regular point of $f\co X\to \R^k$ where $E$ is an extremal subset of $X$.  Suppose $f$ is incomplementable at $p\in E$.
 
 Then there exists a admissible function $h\co X\to\R$ with the following properties
\begin{enumerate}[(i)]
\item  $h(p)=0$.
\item $h$ is strictly concave on $B(p,R)$ for some $R>0$.
\item There are $r>0, A>0$ such that  $h< A$ on $f^{-1}\left(\bar{I}^k(f(p),r) \right)$ and   $f^{-1}\left(\bar{I}^k(f(p),r) \right)\cap  \{h\ge -A\}$ is compact in $B(p,R)$.

\item $h$ has a unique maximum in $B(p,R)\cap f^{-1}(v)$ for all $v\in \bar{I}^k(f(p),r) $.  Let $S$ denote the set of such maximum points.
\item $(f,h)$ is regular on $f^{-1}\left(\bar{I}^k(f(p),r) \right) \cap B(p,R)\backslash S$.
\item $S\subset E$.
\end{enumerate}
\end{lem}
\begin{proof}
The proof is identical to the proof of Lemma~\ref{incompl} except for part (vi) which is new. 
Let $x$ be a  point of max of  $h$ on $E\cap f^{-1}(v)$.  Let $z$ be the unique point of maximum of  $h$ on $ f^{-1}(v)$. If $x\ne z$ then $(f,h)$ is regular at $x$ by Lemma~\ref{incompl}. However, by Lemma~\ref{per-lem-extr}, a point on $E$ is regular for $F\co X\to\R^m$ iff it's regular for $F|_E$. Thus, $(f,h)|_E$ is regular at $x$.  Therefore,  $(f,h)|_E$ is co-lipschitz near $x$ and hence $x$ is not a point of maximum of $h$ on $E\cap f^{-1}(v)$.
\end{proof}

This Lemma easily implies that if $f\co X\to\R^k$ is regular at $p\in E$ where $E\subset X$ is extremal then the local dimension of $E$ near $p$ is $\ge k$ and the equality is  only possible if $F|_E$ is locally bi-Lipschitz near $p$ (cf.  Corollary~\ref {bilip}).

Lemma~\ref{relincompl} also yields the relative local  fibration theorem in exactly the same way as Lemma~\ref{incompl} yields the absolute local fibration theorem.

\begin{proof}[Proof of Theorem~\ref{locrelfib}]
Since the proof is almost identical to the proof of the fibration theorem in ~\cite{Per-Morse} we only give a sketch.

We argue by reverse induction in $k$. Since the base of induction is clear we only have to consider the induction step from $k+1\le n$ to $k$. Let $f\co X\to R^k$ be regular at $p\in E$ where $E\subset X$ is extremal.
If $f$ is complementable at $p$ the statement follows by induction assumption. Suppose $f$ is incomplementable at $p$. Let $h$ be the function provided by Lemma~\ref{relincompl}. 
Suppose for simplicity that $h$ is identically zero on $S$. 
Let $U= f^{-1}\left(\bar{I}^k(f(p),r) \right)\cap\{-A<h<0\}\cap B(p,R)$ and $W=f^{-1}\left(\bar{I}^k(f(p),r) \right)\cap\{-A<h\le 0\}\cap B(p,R)$

Then $U=W\backslash S$ and $(f,h)$ is regular on $U$. Therefore the relative local fibration theorem holds for $(f,h)$ on $U$ by induction assumption.  

By Lemma~\ref{relincompl}, $(f,h)$ is proper on $U$ and hence, 
$(f,h)\co (U,U\cap E)\to \bar{I}^k(f(p),r)\times (-A,0)$ is a relative bundle map by~\cite[6.10]{Sieb}. This means that
$(U,U\cap E)$ is homeomorphic to $(F,B)\times  \bar{I}^k(f(p),r)\times (-A,0)$ respecting $(f,h)$ where $F$ is an MCS-space of $\dim=n-k-1$ and $B$ is a stratified subspace in $F$. By Lemma~\ref{relincompl}, we can extend this homeomorphism to
a homeomorphism $(W, W\cap U)\to (KF,KB)\times \bar{I}^k(f(p),r)$ which proves the induction step.

The general case when $h$ is not  constant on $S$ is handled in exactly the same way as in ~\cite{Per-Morse} and the proof of Key Lemma~\ref{keylemma}  by constructing an auxiliary function $\th$ obtained by shifting $h$ by constants on each of the fiber of $f$ to make it identically zero on $S$.

\end{proof}

Before we can start with the proof of  the Relative Stability Theorem  we first need to observe that by~\cite{Sieb}  the corresponding version of Theorem~\ref{local bundle} and Strong Gluing Theorem~\ref{stronggluing} hold in relative category for pairs of MCS spaces and their stratified subspaces.

The relative version of Theorem~\ref{local bundle} follows from ~\cite[Complement 6.10 to  Union Lemma 6.9]{Sieb} by the same argument as in the proof of the absolute version of Theorem~\ref{local bundle} given by \cite[Corollary 6.14]{Sieb}.

The relative version of the Strong Gluing Theorem still follows from the same deformation of homeomorphism result \cite[Theorem 6.1]{Sieb} which also  covers  relative homeomorphisms.  
Here we'll need the following definition
\begin{defn}
A pair of metric spaces $(X,E)$  is called $\varkappa$-connected if  both $X$ and $E$ (taken with the restricted ambient metric) are $\varkappa$-connected.
\end{defn}

Let us state the Relative Gluing Theorem. For simplicity we only state the unparameterized version.

\begin{relgluingthm}\label{relgluingthm}
Let $(X,E)$ be a stratified pair , $\{U_\a\}_{\a\in \frak A}$ be a finite covering of $X$. Given a function $\varkappa_0$, there exists $\varkappa=\varkappa((X,E), \{U_\a\}_{\a\in \frak A},\varkappa_0)$ such that the following holds:

Given a $\varkappa_0$-connected stratified pair $(\tilde{X},\tilde{E})$,  an open cover of $\tilde{X}$  $\{\tilde{U}_\a\}_{\a\in \frak A}$, a $\delta$-Hausdorff approximation $\phi\co (X,E)\to (\tilde{X},\tilde{E})$ and a family of homeomorphisms $\phi_\a\co 
(U_\a,U_\a\cap E)\to  (\tilde{U}_\a,\tilde{U}_\a\cap \tilde{E})$, $\delta$-close to $\phi$, 

then there exists a homeomorphism $\bar{\phi}\co (X,E)\to  (\tilde{X},\tilde{E})$, $\varkappa(\d)$-close to $\phi$.
\end{relgluingthm}

Observe that under the assumptions of the Relative Stability Theorem, all the elements of the sequence $(X_i,E_i)$  and the limit $(X,E)$ are $\varkappa$-connected by Lemma~\ref{uniform-bi-Lip}.

Furthermore, Lemma~\ref{reg-inner}  and Remark~\ref{cont-inner}  hold for regular level sets of admissible functions on extremal subsets. 
The proof is exactly the same as the proof of  Lemma~\ref{reg-inner} modulo
 Lemma~\ref{uniform-bi-Lip} and the fact that extremal subsets are invariant under all gradient flows.

We are now ready to prove the Relative Stability Theorem.
\begin{proof}[Proof of Theorem~\ref{stab-extr}]

The proof  proceeds by reverse induction on the framing and is, in fact, {\it  exactly} the same as the proof of the usual stability theorem except we make all the arguments relative. Everywhere in the proof substitute $U$ (with various subindices) by $(U,U\cap E)$.  In the proof of the relative version of  Key Lemma~\ref{keylemma}, use Lemma~\ref{relincompl} instead of Lemma~\ref{incompl} whenever necessary. Also, use the Relative Local  Fibration Theorem instead of Local Fibration Theorem~\ref{local bundle} and the Relative Gluing Theorem instead of the Gluing Theorem whenever called for.

\end{proof}

\begin{rmk}
The relative stability theorem trivially implies the following hitherto unobserved fact.  Under the assumptions of the Relative Stability Theorem, $E_i\conv E$ without collapse. Then $\dim E =\dim E_i$ for all large $i$.

A fairly simple direct proof of this statement can be given without using the Relative Stability Theorem. However, we chose not to present it here because this indeed obviously follows from the Relative Stability Theorem.

\end{rmk}

\small
\bibliographystyle{alpha}
\def\cprime{$'$}

\end{document}